\documentclass[11pt]{amsart}

\allowdisplaybreaks[1]
\usepackage{enumerate}

\usepackage{amsmath}
\usepackage{amssymb}
\usepackage{amsfonts}
\usepackage{lmodern}
\usepackage{verbatim,xspace}
\usepackage[usenames]{color}
\usepackage{tikz-cd}
\usepackage{fullpage}
\usepackage[T1]{fontenc}

\usepackage[english]{babel}

\makeindex

\definecolor{cobalt}{rgb}{0.0, 0.28, 0.67}

\definecolor{darkblue}{rgb}{0.0, 0.0, 0.55}

\usepackage[colorlinks,linkcolor=carnelian,citecolor=darkblue,urlcolor=black,hypertexnames=true]{hyperref}

\makeatletter

\def\moverlay{\mathpalette\mov@rlay}

\def\mov@rlay#1#2{\leavevmode\vtop{

                \baselineskip\z@skip \lineskiplimit-\maxdimen

                \ialign{\hfil$#1##$\hfil\cr#2\crcr}}}

 \newtheorem{theorem}{Theorem}[section]

 \newtheorem{corollary}[theorem]{Corollary}

 \newtheorem{lemma}[theorem]{Lemma}

 \newtheorem{proposition}[theorem]{Proposition}

\theoremstyle{definition}

\newtheorem{definition}[theorem]{Definition}

\newtheorem{example}[theorem]{Example}

\theoremstyle{remark}

\newtheorem{remark}[theorem]{Remark}

\newtheorem{fact*}{Fact}

\DeclareMathOperator{\GL}{GL}

\DeclareMathOperator{\SL}{SL}

\DeclareMathOperator{\M}{M}

\DeclareMathOperator{\rest}{Res}
\DeclareMathOperator{\infl}{Inf}

\DeclareMathOperator{\spn}{span}

\DeclareMathOperator{\supp}{supp}

\newcommand{\N}{\mathbb{N}}

\newcommand{\Q}{\mathbb{Q}}
\newcommand{\C}{\mathbb{C}}
\newcommand{\Z}{\mathbb{Z}}

\newcommand{\FF}{\mathbb{F}}

\newcommand{\R}{\mathbb{R}}

\newcommand{\inv}{^{-1}}

\DeclareMathOperator{\Irr}{Irr}
\DeclareMathOperator{\IRR}{IRR}
\DeclareMathOperator{\Ind}{Ind}
\DeclareMathOperator{\Lin}{Lin}
\DeclareMathOperator{\Hom}{Hom}
\DeclareMathOperator{\Aut}{Aut}

\newcommand\al{\alpha}

\newcommand\la{\lambda}

\newcommand{\ls}[1]{{}^{#1}}

\newcommand\bbm{\begin{bmatrix}}
\newcommand\ebm{\end{bmatrix}}

\newcommand{\bpm}{\left( \begin{smallmatrix}}
\newcommand{\epm}{\end{smallmatrix} \right)}

\numberwithin{equation}{section}

\newcommand{\plangle}{\moverlay{(\cr<}}
\newcommand{\prangle}{\moverlay{)\cr>}}

\newcommand{\fpolys}[2]{\mathbb{F}{<}{#1}_1,\ldots, {#1}_{#2}{>}}

\newcommand{\unram}{unramified\xspace}

\newcommand{\semunram}{pseudo-unramified\xspace}

\newcommand{\frats}[2]{\mathbb{C}\plangle{#1}_1,\ldots, {#1}_{#2}\prangle}

\newcommand{\fratsm}[1]{\mathbb{F}\plangle #1\prangle}

\newcommand{\nfrats}[3]{#3\plangle{#1}_1,\ldots, {#1}_{#2}\prangle}

\newcommand{\nfratsm}[2]{#1\plangle #2\prangle}

\newcommand{\mnser}[1]{\mathbb{F}(\!( #1 )\!)}
\newcommand{\mnserdr}[2]{#2(\!( #1 )\!)}

\title{Finite solvable groups with a  rational skew-field of noncommutative real rational invariants}

\author[G. Podlogar]{Gregor Podlogar${}^1$}
\address{
	Institute of Mathematics, Physics and Mechanics,  Ljubljana,
	Slovenia}
	\address{University of Ljubljana, Faculty of Mathematics and Physics Ljubljana, Slovenia}
\email{gregor.podlogar@imfm.si, gregor.podlogar@fmf.uni-lj.si}
\thanks{${}^1$Research done as a part of Doctoral programme at University of Ljubljana, Faculty of Mathematics and Physics, Ljubljana, Slovenia under supervision of Igor Klep}
\date{\today}

\keywords{noncommutative rational invariant, noncommutative Noether's problem, totally \unram groups, multiplicity free restrictions, Clifford theory}

 \subjclass[2020]{Primary 16W22, 20C15; Secondary 16K40, 20C25, 20F22}

\definecolor{carnelian}{rgb}{0.7, 0.11, 0.11}
\definecolor{coolblack}{rgb}{0.0, 0.18, 0.39}
\begin{document}

\begin{abstract}
	For  abelian groups the invariant skew-fields are always rational. 
	We show that for a solvable group  the invariant skew-field is finitely generated. 
	The skew-field invariant under a linear action of a solvable group is rational if the action is well-behaved -- given by a complete representation. 
	We determine the groups that admit such representations and call them totally \semunram.
	In the second part we study  the reach  of  totally \semunram groups and classify  such $p$-groups of rank at most $5$.
\end{abstract}

\maketitle

\tableofcontents

\section{Introduction}
Invariants of  group actions are an important topic that appears in many mathematical areas and beyond mathematics in physics and chemistry.
Classical invariant theory \cite{pop94,spr06} considers polynomials that are preserved under the  action of a  group $G$  
 given by a  homomorphism $G\to \Aut_\FF(\FF[x_1,\dots,x_n])$. 
 Closely related   topic are rational invariants, where one considers rational functions that are invariant under the    action of a  group $G$  
 given by a  homomorphism $G\to \Aut_\FF(\FF(x_1,\dots,x_n))$. 
The  \emph{Noether's problem} asks whether  the subfield of invariants $\FF(x_1,\dots,x_n)^G$ is rational, i.e., purely transcendental over $\FF$. 
The problem depends both on the group $G$ and the base field $\FF$. 
The topic is widely studied \cite{swan83,sal84,C-TS07}.  
Noether's problem has a positive answer in one variable over any field (Lüroth's theorem), in two variables over $\FF=\C$ and for linear actions of abelian groups over $\C$ \cite{end73}. 
There are groups with a negative answer; over $\Q$ even some cyclic groups  (the smallest such is $\Z_{47}$ \cite{len80}), there are examples over $\C$ as well (\cite{sal84}).

We consider a noncommutative version of \emph{Noether's problem}. 
 We replace commutative rational functions with the \emph{free skew-field} $\nfrats{x}{n}{\FF}$ also called the skew-field of \emph{noncommutative rational functions} (in $n$ variables). 
The free skew-field is the universal skew-field of fractions of the free associative algebra of noncommutative polynomials (see \cite{ami66,Coh95}).
The noncommutative Noether's problem then considers the rationality of the skew-field of invariants of a finite group $G$, 
i.e., 
whether the skew field of invariants $\nfrats{x}{n}{\FF}^G$ is isomorphic to the skew-field of noncommutative rational functions 
$\nfrats{y}{m}{\FF}$ for some $m\in \N$.

In \cite{kl20} the authors consider noncommutative rational functions invariant  under a faithful action of finite abelian group  given by linear transformations of variables.
In this case the skew-field of invariants $\frats{x}{n}^A$  is rational in $|A|(n-1)+1$ variables (\cite[4.1]{kl20}). 
For example,
$$
\nfratsm{\C}{x,y}^{S_2}  =\nfratsm{\C}{x+y, (x-y)^2, (x-y)(x+y)(x-y)}\cong  \nfratsm{\C}{y_1,y_2,y_3}.
$$
In fact a closer inspection of the above example shows that we can replace $\C$ by $\Q$ or any other field of characteristic not equal to $2$.

One would like to use the cited theorem recursively on solvable groups; 
given an abelian normal subgroup $N$ of $G$ first compute the $N$-invariants  and proceed with the action of the quotient group $G/N$; however, the subsequent action might not be linear. 
In some cases the idea still works, for example for the action given by a \emph{complete representation} as defined in \cite{kl20} (or see Definition \ref{def-comp}).
If the action of a finite group $G$ is given via a complete representation on the linear span of variables $x_1,\dots, x_n$, then the skew-field of invariants
$\frats{x}{n}^G$  is rational in $|G|(n-1)+1$ variables (\cite[5.1]{kl20}). 
In fact, the proof of the cited theorem serves as an algorithm for expressing free generators of the skew-field of invariants in terms of the initial variables. 
Furthermore,  \cite{kl20} provides a class of groups called \emph{totally \unram} such that their regular representation is complete, 
thus the noncommutative Noether's problem has a (partial) positive answer for them.
Among  totally \unram groups are the symmetric groups $S_3$ and $S_4$. 
Understanding their skew-fields of invariant noncommutative rational functions could shed some light on the theory of noncommutative symmetric rational functions which is, contrary to the commutative case, still far from complete (\cite{gel95}). 

The recursive method  using normal abelian subgroups does not fail completely even for a general finite solvable group $G$; using it we can show that $\frats{x}{n}^G$ is finitely generated as a skew-field over $\C$ (\cite[1.1]{kl20}). 
In contrast,  the ring of noncommutative polynomials invariant under a linear action of a finite group $\fpolys{x}{n}^G$ is almost never finitely generated \cite[6.8.4]{Coh06}.

In this paper we study the noncommutative Noether's problem over $\R$ and prove the real versions of \cite[1.1, 4.1, 5.1]{kl20}.
In the second part of the paper we study the reach of the cited results and our real counterparts.
We introduce \emph{totally \semunram groups}, a generalisation of totally \unram groups, and show they are precisely the groups that admit complete representations, hence we  harvest the full potential of \cite[5.1]{kl20}. 
We  also connect totally \unram groups and totally \semunram groups  with established concepts from group theory and  representation theory.
The  tools for the last part come from Clifford theory.
 Finally, we classify totally \unram and totally \semunram $p$-groups of rank at most $5$. 
These groups are significant as there are $p$-groups of rank $5$ whose fields of commutative rational invariants are not rational \cite{chu10,mor12,hos13}. 
Unsurprisingly none of these groups are totally \semunram.

\subsection{Main results and outline}

  Section \ref{sec-pre} explains the notations and contains the main definitions followed by preliminaries  on noncommutative rational functions, Malcev-Neumann series, Clifford theory and projective representations. 
  
  In Section \ref{sec-real} we study real noncommutative rational functions invariant under a group action given by a complete linear representation and obtain the following main results: 
  \begin{enumerate}
  	\item If the action of a group $G$ on $\nfrats{x}{n}{\FF}$ is nontrivial on $\FF$ and trivial on the variables, then $\nfrats{x}{n}{\FF}^G=\nfrats{x}{n}{\FF^G}$ (Proposition \ref{lem-dve}).
  	  	\item The skew-field of  noncommutative real rational functions invariant under a linear action of an abelian group is rational 
  	(Theorem \ref{prop-ab}).
  	\item The skew-field of  noncommutative real rational functions invariant under an action of a group given via a complete representation is rational 
  	(Theorem \ref{thm-comrat}).
  	\item The skew-field of  noncommutative real rational functions invariant under an action of a finite solvable group is finitely generated.  
  	(Theorem \ref{thm-fingen}).
  	\newcounter{enumTemp}
  	\setcounter{enumTemp}{\theenumi}
  \end{enumerate}
In Section \ref{sec-pse} we study complete representations and totally \semunram groups  with the following results: 
\begin{enumerate}
	\setcounter{enumi}{\theenumTemp}
	\item Totally \semunram groups are precisely the groups that admit complete representations (Theorem \ref{thm-comppseu}).
	\item An example of a group that is not totally \semunram, yet it has a rational skew-field of noncommutative rational invariants (Example \ref{ex-inv}).
	\item A semidirect product $A\rtimes G$ of an abelian group $A$ and a totally \semunram group $G$ is totally \semunram (Corollary \ref{cor-semdir}).
	\setcounter{enumTemp}{\theenumi}
\end{enumerate}
In Section \ref{sec-unr} we study totally \unram groups with main results as follows:  
\begin{enumerate}
	\setcounter{enumi}{\theenumTemp}
	\item Totally \unram groups are closed under quotients   (Proposition \ref{prop-quot}).
	\item Metacyclic groups and semidirect products of abelian groups are totally \unram (Corollaries \ref{prop-cikli} and \ref{prop-semd}).
	\item If $G$ and $H$ are isoclinic finite groups and $G$ is  totally \unram, then $H$ is totally \unram (Proposition \ref{prop-piso}).
	\setcounter{enumTemp}{\theenumi}
\end{enumerate}
Finally, in  Section \ref{sec-nilunr} we study nilpotent totally \unram groups and establish: 
\begin{enumerate}
	\setcounter{enumi}{\theenumTemp}
	\item  A nilpotent totally \unram group is metabelian   (Theorem  \ref{prop-com}).
	\item In Subsection \ref{sub-p} we classify totally \unram  and totally \semunram $p$-groups of rank up to $5$.
\end{enumerate} 
Throughout the paper examples are given to demonstrate the strength of our results.

\section{Definitions and preliminaries}\label{sec-pre}
\subsection{Notation}
Throughout the paper we aim to use  standard notation. 
All considered fields  have characteristic zero.
All considered groups  are assumed to be finite unless stated otherwise.
We denote the set of complex irreducible characters of a group $G$ by $\Irr(G)$ and the set of complex linear characters (characters of degree $1$) by $\Lin(G)$. 
The trivial character of $G$ is denoted by $\tau_G$ or simply $\tau$.
Complex class functions of $G$ are endowed with a scalar product,
$$
\langle \varphi,  \psi\rangle = \frac{1}{|G|}\sum_{g\in G} \varphi(g)\overline{\psi(g)}.
$$ 
The irreducible characters are an orthonormal basis of class functions with respect to this scalar product.
We denote the set of irreducible complex linear representations of a group $G$ by $\IRR(G)$.

A complex  character $\chi$ of a group $G$ is \emph{multiplicity free} if $\langle \chi, \mu\rangle\leq 1$ for every irreducible character $\mu \in\Irr(G)$. 
A complex representation   $\rho$ of a group $G$ is \emph{multiplicity free}  if its character $\chi_\rho$ is multiplicity free. 
Equivalently, $\rho$ is multiplicity free if it is equivalent to a direct sum of pairwise non-equivalent irreducible representations.

For a subgroup $N$ of $G$ and an irreducible character $\mu \in \Irr(N)$ we denote 
the \emph{irreducible characters of $G$ lying over $\mu$} by
$$
\Irr_\mu(G)
=
\{\chi \in \Irr(G) \mid\langle \chi|_N, \mu \rangle > 0\} 
=
\{\chi \in \Irr(G) \mid \langle \chi, \Ind^G_N\mu \rangle>0\}.
$$
By $\IRR_\mu(G)$ we denote the set of irreducible linear representations $\rho$ of $G$ such that their characters satisfy $\chi_\rho\in \Irr_\mu(G)$.

Let $N$ be a normal subgroup of $G$, then
$G$ acts on the characters of $N$. The left action  on a character $\mu$ of $N$ is  defined by $\ls{g}\mu(n)=\mu(g\inv n g)$. 
If $\mu$ is irreducible, then so is $\ls{g}\mu$.

\subsection{Definitions}

We proceed with our main definitions. The definitions of complete representations and totally \unram groups were introduced in \cite{kl20}. 
The original definition of complete representation  is missing the base cases of recursion; we correct this oversight here. 
 Also the definition of $Q\pi$ is slightly changed -- in \cite{kl20} the summand $\pi_{B}\otimes \pi_B$ appears twice. 
Using Proposition \ref{prop-subr} it is easy to see that the definitions are equivalent.

\begin{definition}\label{def-comp}
	A complex linear representation $\pi$ of a finite group $G$ is \emph{complete} if it decomposes as $\pi  =\pi_B \oplus \pi_J$ and 
	there is a nontrivial abelian normal subgroup $N\subseteq G$ such that:
	
	\begin{enumerate}
		
		\item
		
		$\pi_B|_N$ contains exactly the nontrivial linear representations of $N$ as direct summands with multiplicity $1$;
		
		\item $G=N$  or the representation  
		\[
		Q\pi
		=
		\left[ \pi \oplus (\pi_{B}\otimes \pi_B)
		\oplus(\pi_{B}\otimes \pi_J) \oplus (\pi_J \otimes \pi_B)
		\oplus (\pi_{B}\otimes \pi \otimes \pi_B)
		\right]_{N^\tau}
		\]
		is a complete representation of $G/N$.
		Here, for a representation $\rho$,
		$\rho_{N^\tau}$ denotes the summands of $\rho$
		which are trivial on $N$ and thus naturally gives rise to a representation
		of $G/N$.
		\end{enumerate}
	
	A real linear representation of a finite group is \emph{complete} if its complexification is complete.
\end{definition}

\begin{remark}\label{rem-mul}
	Let $\rho\colon G \to \GL(V)$ be a complex linear representation and $N\subseteq G$ an abelian normal subgroup such that $\rho|_N=\bigoplus_{i=1}^m\mu_i$ is multiplicity free. 
	Then $V$ decomposes as a direct sum  $\bigoplus_{i=1}^m V_{\mu_i}$ of one-dimensional subspaces such that 
	$v\in V$ is in $V_{\mu_i}$
	if and only if 
	$\rho(n)v=\mu_i(n)v$ for every $n\in N$. 
	Pick $v_i \in V_{\mu_i}$   then 
	$$
	\rho(g\inv n g)v_i=\mu_i(g\inv n g) v_i = \ls{g}\mu_i(n)v_i.
	$$
	Rearranging yields
	$$
	\rho(n)\left(\rho(g)v_i\right) = \ls{g}\mu_i(n)\left(\rho(g)v_i\right),
	$$
	which shows that $\rho(g)v_i\in V_{\ls{g}\mu_i}$. 
	Hence the matrix of $\rho(g)$ written in any basis $\{b_j \mid j=1,\dots,m\}$ such that $b_j \in V_{\mu_j}$  has precisely one nonzero entry in each row and each column.
\end{remark}

\begin{definition}

	A finite group $G$ is \emph{\unram} over a nontrivial normal abelian subgroup  $N$ if
	for every complex irreducible linear representation $\rho$ of $G$, the restriction $\rho|_N$ is multiplicity free or trivial. 
	
	The group $G$ is \emph{totally \unram} if it is abelian or there exists a nontrivial abelian normal subgroup $N$
	such that $G$ is \unram over $N$ and $G/N$ is totally \unram.
	
\end{definition}
In the definition of a totally \unram group we can exchange representations for their characters. 
Using  Frobenius reciprocity we observe that
a group $G$ is \unram over $N$ if and only if for every nontrivial irreducible representation (character) 
$\mu$ of $N$ the induced representation (character) $\Ind_N^G(\mu)$ is multiplicity free.
\begin{remark}\label{rem-ser}
	We can interpret the recursive condition in the definition of a totally \unram group as follows: a finite group $G$  is totally \unram if there exists a series of normal subgroups
$$
1= N_0   \subsetneq N_1 \subsetneq \dots\subsetneq N_{n-1} \subsetneq N_n=G
$$
such that for every $j=0,\dots, n-1$ the quotient $N_{j+1}/N_{j}$ is abelian and $G/N_{j}$ is \unram over $N_{j+1}/N_{j}$.
\end{remark}

The notion  of "\unram over" is closely related to the so-called Gelfand triples (see \cite{cec20gelf}).
A triple $(G,H,\rho)$  consisting of a group $G$, subgroup $H$ and an irreducible linear representation $\rho$ of $H$ is a \emph{Gelfand triple} if the induced representation $\Ind_H^G\rho$ is multiplicity free. 
Thus a finite group $G$ is \unram over an abelian normal subgroup $N$ if and only if $(G,N,\mu)$ is a Gelfand triple for every  nontrivial irreducible linear representation $\mu$ of $N$.

The next definition is a generalization of totally \unram groups and is, as we will see, tightly connected to complete representations.

\begin{definition}
	A finite group  $G$ is \emph{\semunram over} a nontrivial abelian normal subgroup $N$ if for every irreducible character $\mu\in\Lin(N)$ there exists $\chi \in \Irr_\mu(G)$ such that $\chi|_N$ is multiplicity free.
	
	A finite group  $G$ is \emph{totally \semunram} if it is abelian or there exists a nontrivial abelian normal subgroup $N$	such that $G$ is \semunram over $N$ and $G/N$ is totally \semunram.
\end{definition}

Again we can interchange  representations with  characters and a suitably modified version of Remark \ref{rem-ser} holds for totally \semunram groups.

Clearly every totally \semunram group is solvable.
If $G$ is \unram over $N$, then clearly $G$ is \semunram over $N$. 
Hence every totally \unram group is  totally \semunram.

\begin{remark}
	The definitions of totally \unram and totally \semunram groups allow for a straightforward checking of the properties using  GAP\cite{GAP4}. 
	We use it to work with examples; in the squeal the group with  group ID $[n,m]$ refers to the group that is in GAP summoned by "SmallGroup$(n,m)$". 
\end{remark}

\subsection{Noncommutative rational functions}

The field of (commutative) rational functions is the field of fractions of (commutative) polynomials. 
The passage from noncommutative polynomials i.e., free associative algebra, to the skew-field of noncommutative rational function (also called free skew-field) is not as straightforward.
We introduce terminology and basic concepts surrounding noncommutative rational functions. 
For a longer exposition we refer to \cite{ami66,Coh95,Coh06,kal12,Vol18}. 

A \emph{noncommutative rational expression} is a syntactically valid combination of elements of the base field $\FF$, variables, operations $+,\cdot,$ inverse and parenthesis, for example:
$$\big(2x_1^3x_2^4x_1^5-((x_1x_2 - x_2x_1)^{-1} + 1)^2\big)^{-1}+x_3.$$ 
Such expressions can be evaluated on tuples of square matrices of equal size with coefficients in $\FF$. 
An expression is \emph{nondegenerate} if it is valid to evaluate it on some tuple of matrices. 
Two nondegenerate expressions are equivalent if they evaluate equally whenever both are defined. 
A \emph{noncommutative rational function} is an equivalence class of  a nondegenerate rational expression;
these functions form the free skew-field $\nfrats{x}{n}{\FF}$. 
 The free skew field $\nfrats{x}{n}{\FF}$  is the universal skew-field of fractions of the free algebra  $\fpolys{x}{n}$. 
 It is universal in the sense that any epimorphism from $\fpolys{x}{n}$ to a skew-field $D$ extends to a specialization from  $\nfrats{x}{n}{\FF}$ to $D$.

Another way of constructing  the free skew-field is as the universal localization of free algebra, i.e., we adjoin entries of inverses of all full matrices over the free algebra.
Any noncommutative rational function $r \in \nfrats{x}{n}{\FF}$  can be  represented by a \emph{linear realization}
$$
r = c^*L^{-1}b
$$
where $b,c\in \mathbb{F}^n$ and 
\[L = A_0 + \sum^d_{i=1} A_ix_i\]
for some matrices $A_i \in \M_n(\mathbb{F}).$

We say that a skew-field is \emph{rational} (or \emph{free}) over $\FF$ if it is isomorphic to the free skew-field $\nfrats{x}{n}{\FF}$ for some $n\in \N$.
Variables $x_1,\dots, x_n$ in $\nfrats{x}{n}{\FF}$ generate a free group $\Gamma$ under multiplication. 
The skew-field $\nfrats{x}{n}{\FF}$ is also the universal field of fractions of the group algebra $\FF\Gamma$, hence we also use notation $\nfrats{x}{n}{\FF}=\nfratsm{\FF}{\Gamma}$.

\subsection{Malcev-Neumann series and Connes operator}
Let $\Gamma$  be    the free group  generated by $X=\{x_1,\dots,x_n\}$.
 A \emph{formal (power) series} on $\Gamma$ with coefficients in $\FF$ is a function $S\colon \Gamma \to \FF$.
We denote the set of formal series by $\FF^\Gamma$. 
Any formal series $S$ can be uniquely presented by $\sum_{\omega\in \Gamma} S(\omega) \omega$.
The \emph{support} of a series $S$ is $\supp{S}=\{\omega \mid S(\omega)\neq 0\}$.

Let $\leq$ be any total order of $\Gamma$  compatible with the group structure. 
For an example of such an ordering we refer to  \cite{ber90, reu99}.
The \emph{Malcev-Neumann}  \emph{series} $\mnser{\Gamma, \leq}$ (with respect to the given order) is the set of series  $S \in \FF^\Gamma$ such that their support  is well ordered. Malcev-Neumann series form a skew-field under the pointwise addition and Cauchy product:
$$
ab=\sum_{\omega\in \Gamma}
\sum_{\substack{\alpha,\beta\in\Gamma\\\alpha\beta=\omega}}
a_\alpha b_\beta \omega.
$$
For the sake of brevity we fix the ordering of $\Gamma$ and denote 
$\mnser{\Gamma}=\mnser{\Gamma, \leq}$. 
For more on Malcev-Neumann series we refer to \cite{reu99}.

The rational closure of the free algebra $\fpolys{x}{n}$ or the group algebra $\FF\Gamma$ in $\mnser{\Gamma}$  is isomorphic to $\nfratsm{\FF}{\Gamma}$ regardless of the ordering \cite{lew74,reu99}. 
We say that a series $r\in \mnser{\Gamma}$ is \emph{rational} (over $\FF$) if it belongs to the rational closure of $\fpolys{x}{n}$, i.e., the smallest subring of $\mnser{\Gamma}$ that contains $\fpolys{x}{n}$ and is closed under taking inverse.

Let $M$ be the free monoid generated by $X$. 
Any rational function  $r\in \nfrats{x}{n}{\FF}$
that is defined at $0$ can be expanded  to a series $r \in \FF^M$. 
Conversely, a series  $r \in \FF^M$ represents a rational function if and only if its \emph{Hankel matrix} has  finite rank \cite{ber11,Vol18}.
Rationality in Malcev-Neumann series is a bit more intricate.
Let $\mathcal{G}=\mathrm{Cay}(\Gamma,X)$ be the Cayley graph of $\Gamma$. 
Given a series  $a\in \mnser{\Gamma}$
we have the \emph{Connes operator}  $[\mathfrak{F},a]\colon \FF\mathcal{G} \to \FF^{\mathcal{G}}.$ 
For the definition of the Connes operator we refer to \cite{con94,duc97,lau13}.
A Malcev-Neumann series $a$ is rational if and only if $[\mathfrak{F},a]$ has  finite rank (\cite[12]{duc97}, \cite[2.6]{lau13}).
Given a subfield $\mathbb{K}\subset\FF$ and a series $a \in \mnserdr{\Gamma}{\mathbb{K}}$ the Connes operator 
$[\mathfrak{F},a]\colon \mathbb{K}\mathcal{G} \to {\mathbb{K}}^{\mathcal{G}}$ is equal to the restriction of the Connes operator over $\FF$.

\subsection{Complex noncommutative rational invariants} 
We summarise the results and techniques from \cite{kl20}.
We say that a finite group $G$ acts \emph{linearly} on $\nfrats{x}{n}{\FF}$ if the action is defined by a linear representation $G\to \GL(V)$ where $V=\spn_{\FF}{\{x_1,\dots,x_n\}}$. 
We say that a linear action is \emph{diagonal} if each variable $x_i$ spans an invariant subspace of $V$ i.e., 
$g\cdot x_i = \chi_i(g)x_i$ where $\chi_i$ is a linear character of $G$. 
If $G$ acts  faithfully diagonally on $\nfrats{x}{n}{\FF}$ then $G$ is abelian and $\FF$ is a splitting field of $G$.
Every  linear representation of an abelian group $A$ over a splitting field $\FF$ is equivalent to a direct sum of representation of degree one, thus we can pass from a linear action
of $A$ on  $\nfrats{x}{n}{\FF}$ to a diagonal action via a linear transformation of variables.

Given a faithful diagonal action of a finite abelian group $A$ on  $\nfratsm{\FF}{\Gamma}$ we have a surjective group homomorphism $\Gamma \to A^*(=\Hom(A,\FF^*)\cong A)$ defined by 
$x_i \mapsto \chi_i$. We denote the kernel of this homomorphism by $\Gamma^A$. 
By the Nielsen–Schreier formula, $\Gamma^A$ is a free group of rank $|A|(n-1)+1$. 
\begin{theorem}\label{thm-star}{\cite[4.1]{kl20}}
	If a finite abelian group $A$ acts faithfully diagonally on  $\nfratsm{\FF}{\Gamma}$ then 
	$
	\nfratsm{\FF}{\Gamma}^A=\nfratsm{\FF}{\Gamma^A}.
	$
\end{theorem}
The original statement of the theorem requires $\FF$ to be algebraically closed but allows linear actions,  yet algebraically closed field is only needed to pass from a linear action to a diagonal one.

We continue with invariants of complete representations.
	If the linear action of a group $G$ on $\frats{x}{n}$ is given via a complete representation $\pi=\pi_B\oplus\pi_J$ and $N$ is an abelian normal subgroup from the definition, we can find a "good" set of free generators of $N$-invariants. 
	
	Let $\spn_{\C}\{x_1,\dots, x_n\}=V_B \oplus V_J$ be the decomposition with respect to  $\pi=\pi_B\oplus\pi_J$ and let 
	$$
	\{b_\chi \mid  \chi \in \Irr(N)\backslash\{\tau\}\}
	\quad \text{and} \quad
	\{v_k\mid k=1,\dots,\deg\pi_J\} 
	$$ 
	be  bases of $V_B$ and $V_J$, respectfully,  such that  for each $n\in N$ we have
	$\pi(n)b_\chi=\chi(n)b_\chi$ and   $\pi(n)v_k=\mu_k(n)v_k$ for some $\mu_k \in \Irr(N)$. 
	We also set $b_\tau=1$.
	Then the free generators of $N$-invariants are
	$$
	b_{\chi}b_{\mu}b_{(\chi\mu)\inv}, b_{\theta}v_{k}b_{(\theta\mu_k)\inv}
	$$
	where $\chi$ and $\mu$    run through $\Irr(N)\backslash\{\tau\}$, $\theta$ runs through  $\Irr(N)$ and $k=1,\dots, \deg\pi_J$ (\cite[4.2]{kl20}). 
	Using Remark \ref{rem-mul} we show that $G/N$ acts linearly on these generators via the representation $Q\pi$. 
	The item (2) of the definition then allows us to  continue recursively and conclude that the  skew-field of invariants $\frats{x}{n}^G$ is rational (\cite[5.1]{kl20}).

\subsection{Clifford theory} We give a short overview of  Clifford theory.
For a more thorough and general exposition we refer to  \cite{CR81,Isa94, BE98}.

Let $G$ be a finite group and $N$ an (abelian) normal subgroup. For $\mu\in \Irr(N)$ the \emph{inertia subgroup} is
$I_G(\mu)=\{g \in G \mid \ls{g}\mu =\mu\}$. In this subsection we use $H=I_G(\mu)$.
Pick any left transversal 
$\{g_\alpha \mid \alpha \in G/H\}$ of $H$ in $G$.
  Clifford's theorem (\cite[(11.1)]{CR81}, \cite[(6.5)]{Isa94}, \cite[7.3]{BE98}) states that for any $\chi \in \Irr_\mu(G)$ we have (independently of the   transversal)
$$
\chi|_N = e_\chi\sum_{\alpha \in G/H } \ls{g_\alpha}\mu
$$
where $e_\chi \in \N$ is called the \emph{ramification} of $\chi$ over $N$. 
For any $\chi \in \Irr_\mu(G)$ we have a unique $\theta \in \Irr_\mu(H)$  such that
$$
\chi|_{H}=\sum_{\alpha \in G/H} \ls{g_\alpha}\theta
$$
 and  $\theta|_N=e_\chi \mu$ (\cite[(6.11)]{Isa94},\cite[7.6]{BE98}).
By Frobenius reciprocity we get the dual statement:
$$
\Ind_N^G(\mu)=\Ind^G_{H}(\Ind_N^{H}(\mu))=\sum_{\theta \in \Irr_\mu(H)} e_\theta \Ind^G_{H}(\theta),
$$
 induction   $\theta \mapsto \Ind^G_{H}(\theta)$ gives a bijection from $\Irr_\mu(H)$ to $\Irr_\mu(G)$ and the ramification over $N$  $e_\theta=e_{\Ind^G_{H}(\theta)}$ is preserved.
For our purposes we summarize the findings in the following proposition.

\begin{proposition}\label{prop-lin}
(1)	A group	$G$ is \unram over a nontrivial abelian normal subgroup $N$ if and only if for every nontrivial
 $\mu\in \Irr(N)$ the inclusion $\Irr_\mu(I_G(\mu))\subseteq \Lin(I_G(\mu))$ holds.

(2) A group	$G$ is \semunram over a nontrivial abelian normal subgroup $N$ if and only if for every $\mu \in \Irr(N)$ the intersection $\Irr_\mu(I_G(\mu))\cap\Lin(I_G(\mu))$ is non-empty.
\end{proposition}
\begin{proof}
	(1)	A group $G$ is  \unram over $N$ if and only if for every nontrivial $\mu \in \Irr(N)$ and every $\chi \in \Irr_\mu(G)$
	the ramification $e_\chi$ over $N$  is equal to $1$. 
	  The ramification  $e_\chi$ is the same as the ramification  $e_\theta$ of the unique  $\theta \in \Irr_\mu(I_G(\mu))$ with the property $\Ind^G_{I_G(\mu)}(\theta)=\chi$. 
	We get $\theta|_N=\theta(1)\mu=e_\chi\mu$, hence $e_\chi=\theta(1)=1$ if and only if $\theta\in\Lin(I_G(\mu))$. 
	
	(2) A group $G$ is \semunram over $N$ if and only if for every  $\mu \in \Irr(N)$ there exists 
	$\chi\in \Irr_\mu(G)$ with ramification $e_\chi$ over $N$ equal to $1$.  
	From here we reason as in the proof of (1).
\end{proof}

\subsection{Projective representations}
To give a more palpable description of $\Irr_\mu(I_G(\mu))$ we turn to projective representations. For a detailed discussion we refer to \cite[Ch.\! 6]{BE98}.

Let $V$ be a finite dimensional (complex) vector space.
\emph{A (complex) projective representation}  of a group $G$  is a mapping $P\colon G \to \GL(V)$ 
satisfying 
$$
\forall x,y\in G: P(x)P(y)=\pi(x,y)P(xy) 
$$
 for some $\pi\colon G \times G \to \C^*$. 
 We call $\pi$  a \emph{factor set} of $P$ and  $P$  a \emph{$\pi$-representation}. 
The \emph{degree} of representation $P$ is the dimension of $V$.
A $\pi$-representation is \emph{irreducible} if it does not have any nontrivial invariant subspaces. 
We denote the set of  irreducible $\pi$-representations of $G$ by $\IRR^\pi(G)$.

Any 
factor set $\pi$ satisfies  the $2$-cocycle condition:
$$
\forall x,y,z\in G: \pi(x,y) \pi(xy,z)= \pi(x,yz) \pi(y,z).
$$
Conversely any  mapping $\pi\colon G\times G \to \C^*$  satisfying the $2$-cocycle condition ($2$-cocycle) is a factor set of some projective representation.
The factor sets   equipped with  pointwise  multiplication form the abelian  $2$-cocycle group $Z^2(G,\C^*)$   (with trivial action on $\C^*$).
The factor sets $\pi$, $\pi'$ are \emph{associated} if there exists a function $\la \colon G \to \C^*$ such that
$$
\forall x,y\in G: \pi'(x,y)=\frac{\la(x)\la(y)}{\la(xy)}\pi(x,y).
$$
The factor sets that are associated to the trivial factor set form the subgroup $B^2(G,\C^*)$ of  $2$-coboundaries.
The second cohomology group
$$
M(G)=H^2(G,\C^*)=Z^2(G,\C^*)/B^2(G,\C^*)
$$ 
is also called the \emph{Schur multiplier} of $G$. 
We denote the equivalence class of a factor set $\pi \in Z^2(G,\C^*) $ by $[\pi]\in M(G)$.
We  get another equivalent definition of the Schur multiplier using  Hopf's formula $M(G)=H_2(G,\Z)\cong (R\cap [F,F])/[R,F]$, 
where $F$ is a free group and $R$ a normal subgroup such that $G\cong F/R$. 
For a thorough exposition on the Schur multiplier  we refer to \cite{Kar87}.  

The next proposition follows directly from  $M(G)\cong (R\cap [F,F])/[R,F]$.

\begin{proposition}\label{prop-cycshr}
	The Schur multiplier of a finite cyclic group is trivial.
\end{proposition}

Projective representations $P\colon G \to \GL(V)$ and $P'\colon G \to \GL(V')$ are \emph{linearly equivalent} if there exists a linear isomorphism  
$S \colon V \to V'$ such that
$P'(g)=SP(g)S\inv$  for every $g \in G$. 
Projective representations $P$ and $P'$ are \emph{projectively equivalent} if there exists a function $\la\colon G \to \C^*$ such that $\la P$ and $P'$ are linearly equivalent.
We note that projectively equivalent representations have associated factor sets.
 Conversely there is a bijection between $\IRR^\pi(G)$ and $\IRR^{\pi'}(G)$ if $\pi$ and $\pi'$ are associated. 
 Namely, if 
$\pi'(g,h)=\frac{\la(g)\la(h)}{\la(gh)}\pi(g,h)$, then
the bijection is given by $P\mapsto \la P$.

Some properties  of $\IRR^\pi(G)$ are determined by the class $[\pi]\in M(G)$.
One such property is described in the next lemma.
\begin{lemma}\label{lem-one}
A group $G$ has a $\pi$-representation of degree one if and only if $[\pi]=1\in M(G)$.
\end{lemma}
\begin{proof}
If $P$ is a $\pi$-representation of degree one,  we get 
$
\pi(g,h)=P(g)P(h)/P(gh)
.
$
Conversely if  $
\pi(g,h)=\la(g)\la(h)/\la(gh),
$
then $\la$ is a $\pi$-representation of degree one.
\end{proof}

We now describe  a connection between $\IRR_\mu(I_G(\mu))$ for an irreducible character $\mu \in \Irr(N)$ of an abelian normal subgroup $N$ of $G$ and some projective representations of $I_G(\mu)/N$.
In the sequel we denote $H=I_G(\mu)$.

Let $\rho\colon H \to \GL(V)$ be a linear representation of $H$ of degree $n$ with character $\chi_\rho \in \Irr_\mu(H)$ for some $\mu \in \Irr(N)$. 
Further we choose a left transversal $\{h_\alpha\mid \alpha \in H/N \}$ of $N$ in $H$. 
The transversal defines a 2-cocycle $f\in Z^2(H/N,N)$ by 
$h_\alpha h_\beta=f(\alpha,\beta)h_{\alpha\beta}$. 
We define a mapping $\Gamma\colon H/N \to \GL(V)$ by 
$$
P(\alpha)=\rho(h_\alpha).
$$
Then $P$ is an irreducible projective representation of $H/N$ of degree $n$ with the factor set $\pi(\alpha, \beta)=\mu(f(\alpha,\beta))$. 
We say that $P$ is a \emph{descent} of $\rho$ to $H/N$.
A different choice of transversal defines an associated factor set and  a projectively equivalent representation.
 
We can reverse the process. 
Let $P\colon H/N \to \GL(V)$ be an irreducible projective representation of $H/N$ with a  factor set $\pi(\alpha, \beta)=\mu(f(\alpha,\beta))$.
We define a mapping $P\colon H \to\GL(V)$ by 
$$
\rho(nh_\alpha)=\mu(n)P(\alpha).
$$
Then $\rho$ is an irreducible linear representation of $H$ of degree $n$ and its character lies in $\Irr_\mu(H)$. 
We say that $\rho$ is a \emph{lift} of $P$ to $H$.

We summarise the above discussion in a lemma.
\begin{lemma}\label{lem-deg}
	Let $N$ be an abelian normal subgroup of $G$ and $\mu\in \Irr(N)$ an irreducible character. 
	Further let $f\in  Z^2(I_G(\mu)/N,N)$ be any $2$-cocycle defining the extension from $I_G(\mu)/N$ to $I_G(\mu)$.
	Then there is a degree preserving bijection between $\IRR_\mu(I_G(\mu))$  and $\IRR^\pi(I_G(\mu)/N)$, where $\pi=\mu\circ f\in Z^2(I_G(\mu)/N,\C^*)$.
\end{lemma}

\subsection{Inflation-restriction exact sequence}

Motivated by Lemma \ref{lem-deg} we investigate the mapping $\mu \mapsto [\mu\circ f]$. 
This mapping appears in  a case of the inflation-restriction exact sequence. 
The discussed sequence  is thoroughly examined in \cite[Ch.6.§5]{BE98}.

Let $N$ be an (abelian) normal subgroup of $H$.
We denote by
$$
\Lin^H(N)=\{\mu \in \Lin(N) \mid \mu^h=\mu \text{ for any } h \in H\}=\{\mu \in \Lin(N) \mid I_H(\mu)=H \}
$$ 
\emph{the linear characters of $N$ invariant under the action of $H$}. 
There is an exact sequence
\begin{equation}\label{ex-lhs}
\begin{tikzcd}
		1 \ar[r]&\ar[r,"\infl"]	\Lin(H/N)  &\ar[r,"\rest"] \Lin(H)&\Lin^H(N)\ar[r,"T_H"]&M(H/N)
\end{tikzcd}
\end{equation}
 where  $T_H$ is given by $T_H(\mu)=[\mu\circ f]$ with $f\in Z^2(H/N,N)$ being any $2$-cocycle that describes the extension of $H/N$ to $H$.
We also note that $\Lin(H)\cong H/H'$ and $\Lin^H(N)\cong N/[N,H]$.

\section{Real noncommutative rational invariants}\label{sec-real}

In this section  we show that the noncommutative real rational functions invariant under the linear action of an abelian group are rational over $\R$. 
Then we extend the result to actions given by complete representations. 
At the end we show that the skew-field of rational invariants of a finite solvable group is finitely generated.

We derive a technique that will allow us to pass from the invariants over $\C$ to invariants over $\R$ using the action of the Galois group  $\mathrm{Gal}(\C/\R)\cong \Z_2$.
For this we consider group actions on noncommutative rational functions that are nontrivial on the base field and trivial on the variables.
 We use   Malcev-Neumann series to study   functions invariant under such actions.

\begin{proposition}\label{lem-dve}\begin{enumerate}[(1)]
\item Let the action  of a (possibly infinite) group  $G$ on $\mnser{\Gamma}$ be given by an action on $\FF$ and trivial action on $\Gamma$, then $\mnser{\Gamma}^G=\mnserdr{\Gamma}{\FF^G}$.

\item Let $\mathbb{K}\subset \FF$ be fields, then 
$\mnserdr{\Gamma}{\mathbb{K}}\cap\nfratsm{\FF}{\Gamma}=\nfratsm{\mathbb{K}}{\Gamma}$.

\item Let  the action of a (possibly infinite) group $G$ on $\nfratsm{\FF}{\Gamma}$ be given by an action on  $\FF$  and a trivial action  on  $\Gamma$, 
then $\nfratsm{\FF}{\Gamma}^G=\nfratsm{\FF^G}{\Gamma}.$
\end{enumerate}
\end{proposition}
\begin{proof}
(1) If a series $\sum_{\omega\in \Gamma}a_{\omega}\omega \in \mnser{\Gamma}$ is invariant we get  
$$
g\cdot\sum_{\omega\in \Gamma}a_{\omega}\omega
=\sum_{\omega\in \Gamma}(g\cdot a_{\omega})\omega
=\sum_{\omega\in \Gamma}a_{\omega}\omega
$$
for every $g\in G$, hence the coefficients are invariant.

(2) Take a series  $a\in \mnserdr{\Gamma}{\mathbb{K}}\cap\nfratsm{\FF}{\Gamma}$. 
The rank of Connes operator $[\mathfrak{F},a]$ over $\FF$ is finite, hence the rank of Connes operator over $\mathbb{K}$ is finite as well, thus $a$ is rational over $\mathbb{K}$.

(3)	We embed $\nfratsm{\FF}{\Gamma}$ into the skew-field of Malcev-Neumann series $\mnser{\Gamma}$. 
The action of $G$ on $\nfrats{x}{n}{\FF}$ extends to an action on $\FF(\!(\Gamma)\!)$ given by
$$g\cdot\sum_{\omega\in \Gamma}a_{\omega}\omega=\sum_{\omega\in \Gamma}(g\cdot a_{\omega})\omega.$$ 
Now we apply (1) and (2).
\end{proof}

Given a rational expression or a realization of a noncommutative complex rational function   it is not immediately obvious that the function can be written as a sum of real and imaginary part. 
\begin{corollary}
	Any  $r\in \nfrats{x}{n}{\C}$ can be written as $r_R+ir_I$ with $r_R,r_I\in \nfrats{x}{n}{\R}$.
\end{corollary}
\begin{proof}
	Let $\Z_2$ act on $\C$ by the complex conjugation and trivially on the variables. 
	We denote the action of the nontrivial element of $\Z_2$ on the function $r$  by $\overline{r}$.
	We define $r_R=\frac{1}{2}(r+\overline{r})$ and $r_I=-i\frac{1}{2}(r-\overline{r})$. 
	By Proposition \ref{lem-dve}, both functions are in $\nfrats{x}{n}{\R}$.
\end{proof}
We can always translate certain actions of $\Z_2$ on $\frats{x}{n}$ to  actions that fit the premise of Proposition \ref{lem-dve}.
\begin{lemma}\label{lem-z2}
	Let the action of $\Z_2$ on $\fratsm{\Gamma}$ be defined by a group automorphism of $\Gamma$ and a field automorphism of $\FF$. 
	There exist free generators of $\fratsm{\Gamma}$ such that $\Z_2$ acts diagonally on them.
	If the action on $\FF$ is nontrivial, then there exist free generators of $\fratsm{\Gamma}$ such that $\Z_2$ acts trivially on them.
\end{lemma}
\begin{proof}
Let $\theta\colon \Gamma \to \Gamma$ be a group automorphism of order  two.  
By \cite{mcc80}, there exist free generators $\{w_1,\dots,w_n\}$ of $\Gamma$ such that $\theta w_j= w_j\inv$ or $\theta w_j= A_j w_j B_j$ where $A_j,B_j$ only depend on $w_k$ for $k<j$ and $\theta A_j = A_j\inv, \theta B_j = B_j\inv$.
We replace the generators $w_j$ for which $\theta w_j= w_j\inv$ with $y_j=(w_j-1)(w_j+1)\inv$ and get $\theta y_j = -y_j$ and we replace generators $w_j$ for which 
$\theta w_j= A_j w_j B_j$ by $y_j=(1+A_j) w_j (1+B_j)$ to get $\theta y_j =y_j$.
	
If the automorphism $\theta$ is nontrivial  on $\FF$ then there exists $c\in \FF$ such that $\theta c =-c$.	
We replace the free generators $y_j$ such that  $\theta y_j=- y_j$ by $y_j'=cy_j$  and get 
$\theta y_j'=\theta c \theta y_j=y_j'$.
\end{proof}
\begin{corollary}\label{cor-rat}
	Let the action of $\Z_2$ on $\fratsm{\Gamma}$ be given via  automorphisms of $\Gamma$ and $\FF$. 
	If the action is trivial on $\FF$, then $\fratsm{\Gamma}^{\Z_2}$ is rational  over $\FF$.
	If the action is nontrivial of $\FF$, then $\fratsm{\Gamma}^{\Z_2}$ is rational  over $\FF^{\Z_2}.$
\end{corollary}
\begin{proof}
	First we apply Lemma \ref{lem-z2}. 
	If the action on  $\FF$ is nontrivial we finish with Lemma \ref{lem-z2}, otherwise with Theorem \ref{thm-star}. 
\end{proof}

The method to determine real noncommutative invariants consists of expanding the action to rational functions over $\C$, computing complex invariants using   \cite[4.1]{kl20} and then pushing the results back to noncommutative rational functions over $\R$ using the above results. 
\begin{theorem}\label{prop-ab}
	Let a finite abelian group  $A$ act faithfully linearly on $\nfrats{x}{n}{\R}$, then
	 $$\nfrats{x}{n}{\R}^{A}\cong \nfrats{w}{|A|(n-1)+1}{\R}.$$
\end{theorem}
\begin{proof}
	We extend the action of $A$ to $\nfrats{x}{n}{\C}$
	and let $\Z_2$ act on $\C$ by the complex conjugation and trivially on the variables.  The actions commute, hence we get an action of the group 
	$A\times \Z_2$ and invariants $\nfrats{x}{n}{\R}^{A}=\nfrats{x}{n}{\C}^{A\times \Z_2}.$
		
	We note that every representation of an abelian group over $\R$ is equivalent to a direct sum of irreducible representations of dimension at most $2$.
	One dimensional representations are of the form $g\cdot x = \pm x$. 
	Two dimensional irreducible representations are equivalent to a representation of the form
	$$
	g\cdot \bbm x \\ y \ebm =
	\bbm \cos \phi & -\sin \phi\\ \sin \phi & \cos \phi \ebm
	\bbm x \\ y \ebm
	$$
	and they are diagonalized over $\C$ in the variables $x+iy$ and $x-iy$, which are permuted by the action of $\Z_2$. 
	Therefore we can choose a set of free generators $X$  of $\nfrats{x}{n}{\C}$  such that $A$ acts diagonally on them and $X$ is closed under the action of $\Z_2$.
	Let
	$\Gamma$ be the group generated by $X$. Then $\Gamma$ is closed under the action of $\Z_2$,  hence $\Gamma^A$ is also closed under the action of $\Z_2$.
	By Corollary \ref{cor-rat} the skew-field of invariants $\nfrats{x}{n}{\R}^{A}=\nfratsm{\C}{\Gamma^A}^{\Z_2}$ is rational over $\R$ in $|A|(n-1)+1$ variables.
\end{proof}

\begin{example}\label{exa-z3}
	Let $\Z_3$ act on  $\nfratsm{\R}{x,y,z}$ via a cyclic permutation of variables ($x\mapsto y \mapsto z \mapsto x$). 
	Take $\omega=-\frac{1}{2}+i\frac{\sqrt{3}}{2}$. 
The action is diagonalised over $\C$ in variables 
	$$a=x+y+z, b=x+\omega y+ \overline\omega z, c=x+\overline\omega y+ \omega z$$
	and the free generators of $\nfratsm{\C}{x,y,z}^{\Z_3}$ are
	$$
	a, bc, bac, cb,cab, b^3,c^3.
	$$
	We use a linear transformation to get the free generators of $\nfratsm{\R}{x,y,z}^{\Z_3}$:
	\begin{align}\label{eq-gen}
		w_1=a, w_2=bc+cb, w_3=i(bc-cb)&, w_4=bac+cab,
		\\ w_5=i(bac-cab),w_6=b^3+c^3&,  w_7=i(b^3-c^3) \nonumber
	\end{align}
To  express the generators in the initial variables $x,y,z$ we first introduce some notation:
			\begin{align*}
		f_1&(A,B,C)=AB+BC+CA, \quad
		f_2(A_1,A_2,A_3)= \sum_{\sigma\in S_3} A_{\sigma(1)}A_{\sigma(2)}A_{\sigma(3)},\\
			f_3&(A,B,C)=A(A+B+C)B+B(A+B+C)C+C(A+B+C)A.
	\end{align*}
The free generators $w_j$ are then expressed as follow:
	\begin{align*}
	w_1={}&x+y+z,\\
	w_2={}&2\left(x^2+y^2+z^2\right)-f_1(x,y,z)-f_1(x,z,y),\\ 
	w_3={}&\sqrt{3}(f_1(x,y,z)-f_1(x,z,y)),\\
w_4={}&3\left( x(x+y+z)x+y(x+y+z)y+z(x+y+z)z\right)-(x+y+z)^3\\
w_5={}&\sqrt{3}\left(f_3(x,y,z)-f_3(x,z,y) \right)\\
	w_6={}&2\left(x^3+y^3+z^3+
f_2(x,y,z)\right)-\\
	&\frac{1}{2}\left(f_2(x,x,y)+f_2(x,x,z)+f_2(x,y,y)+f_2(x,z,z)+f_2(y,y,z)+f_2(y,z,z)\right),\\
	w_7={}&\frac{\sqrt{3}}{2}\left(f_2(x,x,z)+f_2(x,y,y)+f_2(y,z,z)-f_2(x,x,y)-f_2(x,z,z)-f_2(y,y,z)\right).
\end{align*}
\end{example}

Next we consider the invariants of an action given via a complete representation.
\begin{theorem}\label{thm-comrat}
	Let a finite group $G$ act on  $\nfrats{x}{n}{\R}$ via a  complete representation, then 
	$$\nfrats{x}{n}{\R}^{G}\cong \nfrats{w}{|G|(n-1)+1}{\R}.$$
\end{theorem}
\begin{proof}
	Denote the complete representation by $\pi$.
	We extend the action of $G$ to $\nfrats{x}{n}{\C}$ and let $\Z_2$ act on $\C$ by the complex conjugation and trivially on the variables.
	Let $N$ be the normal abelian subgroup of $G$ and $\spn_{\C}\{x_1,\dots, x_n\}=V_B \oplus V_J$ the decomposition with respect to $\pi$ according to the definition of complete representation.
	We get the  basis of $V$
	$$
	\{b_\chi \mid  \chi \in \Irr(N)\backslash\{\tau\}\}
	\cup
	\{v_k\mid k=1,\dots,\deg\pi_J\}.
	$$
	that diagonalize $\pi|_N$.
	We denote $b_\tau = 1$. 
	Since every real irreducible character of an abelian group decomposes over $\C$ as a sum of conjugate linear characters we can choose a basis such
	 that $\overline{b_\chi}=b_{\bar\chi}=b_{\chi\inv}$ and $\overline{v_k}=v_{k'}$ for some $k'$.
	 By  \cite[4.2]{kl20},
	 $$
	 \left\{
	 b_{\chi}b_{\mu}b_{(\chi\mu)\inv}, b_{\theta}v_{k}b_{(\theta\mu_k)\inv} 
	 \mid \chi,\mu \in \Irr(N)\backslash\{\tau\}, \theta \in  \Irr(N), k=1,\dots,\deg\pi_J
	 \right\}
	 $$
	 are the free generators of $\nfrats{x}{n}{\C}^N$ such that $G/N$ acts linearly on them.
	 We use a linear transformation to get   free generators
	 $$
	 b_{\chi}b_{\mu}b_{\overline{\chi\mu}}+
	 b_{\bar\chi}b_{\bar\mu}b_{\chi\mu}, \quad
	  i(b_{\chi}b_{\mu}b_{(\chi\mu)\inv}-
	  b_{\bar\chi}b_{\bar\mu}b_{\chi\mu})
	 $$
	 and
	 $$
	 b_{\theta}v_{k}b_{\overline{\theta\mu_k}}+
	 b_{\bar\theta}\overline{v_{k}}b_{\theta\mu_k}, \quad
	 i(b_{\theta}v_{k}b_{\overline{\theta\mu_k}}-
	 b_{\bar\theta}\overline{v_{k}}b_{\theta\mu_k}),
	 $$
where $\chi$ and $\mu$    run through $\Irr(N)\backslash\{\tau\}$, $\theta$ runs through  $\Irr(N)$ and $k=1,\dots, \deg\pi_J$. 
We note that some expressions appear twice and some are equal to $0$; these we omit.
Then $\Z_2$ acts trivially on these free generators, hence they are also free generators of $\nfrats{x}{n}{\R}^N$ by Proposition \ref{lem-dve}.
The quotient $G/N$ acts linearly on them via a representation that is equivalent to $Q\pi$, thus we can continue with recursion.
\end{proof}

\begin{example}
	Let the symmetric group $S_3$ act on $\nfratsm{\R}{x,y,z}$ by permuting the variables. 
	The invariants of the cyclic normal  subgroup $\langle (123) \rangle\cong \Z_3$  are computed in Example \ref{exa-z3}. 
	The action of the quotient $\Z_2 \cong S_3/\langle (123) \rangle$  is given by the action of the cycle $(12)$, that is defined by
	$$
	a\mapsto a, \ b\mapsto \overline\omega c, \  c \mapsto \omega b,
	$$ 
	hence the action on the free generators given by \eqref{eq-gen}  is 
$$
w_1 \mapsto w_1, w_2 \mapsto w_2, w_3 \mapsto -w_3, w_4 \mapsto w_4, w_5 \mapsto -w_5, w_6 \mapsto w_6, w_7 \mapsto -w_7
$$
and the free generators of
$\nfratsm{\R}{x,y,z}^{S_3}$ are
$$
w_1, w_2 ,w_4, w_6, w_3w_3, w_3w_5, w_3w_7, w_5w_3, w_7w_3, w_3w_1w_3, w_3w_2w_3, w_3w_4w_3, w_3w_6w_3.
$$
With a bit of further effort we can show that these are also the free generators of $\nfratsm{\Q}{x,y,z}^{S_3}$.
\end{example}
 
The computing method fails for a general linear action of a solvable group, however we can still use it to show that the skew-field of noncommutative real rational invariants of a solvable group is finitely generated (as a skew-field over $\R$),
whereas  the ring of noncommutative polynomials invariant under a linear action of a finite group is almost never finitely generated \cite[6.8.4]{Coh06}.  

We begin with abelian groups.
\begin{lemma}\label{lem-fin}
	Let $D$ be a finitely generated skew-field over $\R$ and let a finite abelian group  $A$ act on $D$ via a homomorphism $A\to\Aut_{\R}(D)$. 
	Then $D^A$ is a finitely generated skew-field over $\R$.  
\end{lemma}
\begin{proof}
	Let $d_1,\dots, d_m$ be the generators of $D$.
	We define an action of $A$ on the free skew-field $E=\nfratsm{\R}{x_{ig}\mid i=1,\dots,m, g\in A}$ by $gx_{ih}=x_{i(gh)}$. 
	Then the specialization $\psi \colon E \to D$ defined by $\psi(x_{ig})=gd_i$ satisfies $\psi(gy)=g\psi(y)$ and its domain  is closed under the action of $A$. 
	By Theorem~\ref{prop-ab}, $E^A$ is rational in finitely many variables.
	Further, the action of $A$ on $E$ is given by a direct sum of copies of the regular representation of $A$, hence it contains every representation of $A$. 
	Thus we proceed as in the proof of Theorem \ref{thm-comrat} to find free generators of  $E^A$ that are polynomial in the initial variables, hence they are in the domain of $\psi$.
	 
	We restrict the specialization $\psi$ to $E^A$. 
	Any $d\in D^A$ is the image of some $r\in F$, but then it is also the image of $1/|A|\sum_{g\in A} g r\in E^A$ which shows that $\psi|_{E^A}$ is surjective.
	Hence, $D^A$ is generated by the images of the generators of $E^A$.
\end{proof}
We apply the lemma inductively to deduce the result for solvable groups. 
\begin{theorem}\label{thm-fingen}
		Let $D$ be a finitely generated skew-field over $\R$ and let a finite solvable group  $G$ act on $D$ via a homomorphism $G\to\Aut_{\R}(D)$. 
	Then $D^G$ is a finitely generated skew-field over $\R$.  
	In particular this holds for $D=\nfrats{x}{n}{\R}^{G}$.
\end{theorem}
\begin{proof}
	Let $N$ be a nontrivial abelian subgroup of $G$. 
	By Lemma \ref{lem-fin}, the skew-field  $D^N$ is finitely generated. 
	We continue with the action of $G/N$ on $D^N$  and  conclude the proof by induction.
\end{proof}

\section{Complete representations and totally \semunram groups}\label{sec-pse}

In this section we study complete representations and totally \semunram groups.
Our first nonabelian examples of totally \semunram (and totally \unram) groups are dihedral groups.  
\begin{example}\label{ex-comp}
	The dihedral group $D_{2n}=\langle a,b \mid a^n =  b^2 =abab= e\rangle$ ($n\geq3$) has irreducible representations of degree one and two.
	The representatives of two-dimensional irreducible representations  are given by 
	$$
	\pi_\omega \colon a \mapsto \bbm
	\omega & 0\\
	0 & \omega\inv
	\ebm,
	\ b \mapsto
	\bbm
	0 & 1\\
	1 & 0
	\ebm,
	$$ 
	where $\omega$ is a $n$-th root of unity such that $\omega\neq\omega\inv$. 
	The restriction of $\pi_\omega$ to the normal subgroup $\langle a \rangle$ is equivalent to 
	$\pi_\omega|_{\langle a \rangle} \cong \mu_\omega \oplus \mu_{\omega\inv}$, where $\mu_x(a)=x$. 
	Clearly the restriction is multiplicity free, therefore $D_{2n}$ is \unram  over $\langle a \rangle$. 
	The quotient $D_{2n}/\langle a \rangle$ is abelian, hence, $D_{2n}$ is totally \unram and also totally \semunram.
	If $n=3$, the representation  $\pi_\omega$ is complete, otherwise it is not.
\end{example}
The above example shows that the standard representation of $D_6\cong S_3$ is complete. 
We will show that the same holds for the standard representation of $S_4$. 
First we prove an easy proposition.
\begin{proposition}\label{prop-subr}
 (1) Let $\pi$ be a subrepresentation of $\rho$ and assume $\pi$ is complete, then $\rho$ is complete.
 
 (2) Let $\rho$ be a subrepresentation of $\pi$  and suppose $\pi\oplus\rho$ is complete, then $\pi$ is complete.
\end{proposition}
\begin{proof}
	(1)	Write $\pi=\pi_B\oplus\pi_J$ as guaranteed in the definition of a complete representation, further write  $\rho = \pi \oplus \pi_0$. 
	Then we can decompose $\rho=\rho_B\oplus\rho_J$ where $\rho_B=\pi_B$ and $\rho_J=\pi_J\oplus\pi_0$. 
	Clearly this decomposition satisfies item (1) from the definition. 
	The item (2) is clearly true if we are at the last step of the recursion.
	Otherwise  $Q\pi$ is a subrepresentation of $Q\rho$ and we finish by recursion.
	
	(2) We can assume that $\rho$ is irreducible. 
	Write $\pi\oplus\rho=\pi_B\oplus \pi'_J$. 
	Item (1) of the definition of  a complete representation shows  that $\rho$ appears in $\pi_B$ with multiplicity at most $1$, 
	hence we can write $\pi'_J=\pi_J\oplus \rho$ and $\pi=\pi_B\oplus \pi_J$. Then we compute
	$$
	Q(\pi\oplus \rho)= Q\pi \oplus 
	\left[ \rho \oplus (\pi_B \otimes \rho) \oplus (\rho \otimes \pi_B) \oplus (\pi_B \otimes \rho \otimes \pi_B)
	\right]_{N^\tau},
	$$
	since every summand of 
	$$
	\left[ 
		\rho \oplus (\pi_B \otimes \rho) \oplus (\rho \otimes \pi_B) \oplus (\pi_B \otimes \rho \otimes \pi_B)
	\right]_{N^\tau}
	$$
	is contained in $Q\pi$ as a subrepresentation we can finish with recursion.
\end{proof}
\begin{example}
		The symmetric group $S_4$ contains an abelian normal  subgroup 
		$$V=\{e, (12)(34), (13)(24), (14)(23)\},$$ which is isomorphic to the Klein four-group $\Z_2\times\Z_2$. 
	The character tables  of $S_4$ and $V$ are:
		$$
		\begin{array}{c|ccccc}
			S_4&\{\}&\{2\}&\{2,2\}&\{3\}&\{4\}\\ \hline
			\tau_{S_4}&1&1&1&1&1\\
			\chi_1&1&-1&1&1&-1\\
			\chi_2&3&1&-1&0&-1\\
			\chi_3&3&-1&-1&0&1\\
			\chi_4&2&0&2&-1&0\\
		\end{array}\qquad
		\begin{array}{c|rrrr}
			V&&&& \\ \hline
			\tau_V&1&1&1&1\\
			\mu_1&1&-1&1&-1\\
			\mu_2&1&1&-1&-1\\
			\mu_3&1&-1&-1&1\\
		\end{array}
	$$
The character of the standard representation $\rho$ of $S_4$ is $\chi_2$. 
Its restriction to $V$ is  $\chi_2|_V=\mu_1+\mu_2+\mu_3$, thus $\rho$ satisfies item (1) of the definition of a complete representation. 
It remains to show that $Q\rho$ is complete. 
For this we compute $\chi_2\otimes\chi_2=\tau+\chi_2+\chi_3+\chi_4$. 
Since $\chi_4$ is induced by the character of the standard representation of $S_3\cong S_4/V$, the representation $Q\rho$ contains the standard representation of $S_3$ as a subrepresentation and is therefore complete. 

Tracing through the character table it is  not hard to see that $S_4$ is \unram over $V$, thus $S_4$ is totally \unram.
We cannot extend these examples to $S_5$ as it is not solvable.
\end{example}

As promised in the introduction  we connect  complete representations and totally \semunram groups  in the next theorem.
\begin{theorem} \label{thm-comppseu}
	(1)	The regular representation of  a totally \semunram group is complete.
	
	(2)	If a group admits a complete representation, then it is totally \semunram.
	
	(3)	A group  is totally \semunram if and only if it admits a  complete representation.
\end{theorem}
\begin{proof}
	(1) Let $G$ be totally \semunram and let
	$R_G$ be the regular representation. Let $G$ be \semunram over an abelian normal subgroup $N$.
	Partition $\Irr(N)$ into equivalence classes of the form $[\mu]=\{\ls{g}\mu\mid g \in G\}$. For each such class of nontrivial characters pick   
	$\chi_{[\mu]}\in \Irr(G)$ such that 
	$\langle\chi_{[\mu]}|_N,\theta\rangle=1$ for any (and hence all) $\theta\in[\mu]$ as guaranteed in the definition of "\semunram over". 
	Write $\chi_B$ for the sum of these characters and let $\pi_B$ be a subrepresentation of $R_G$ with character $\chi_B$. 
	Then write $R_G=\pi_B\oplus\pi_J$. 
	By construction, $\pi_B|_N$ contains all nontrivial representations of $N$ with multiplicity one.
	Notice that $\pi_J$ contains a regular representation of $G/N$ as a subrepresentation and therefore the representation $Q\pi$ contains a regular representation of $G/N$. By Proposition \ref{prop-subr} we can proceed with recursion.
	
	(2) Let $\pi$ be a complete representation of $G$ and let $N$ be the abelian normal subgroup from the definition of complete representation. 
	Then the characters of the irreducible summands of $\pi_B$ are the characters required to show that $G$ is \semunram over $N$. 
	We conclude the proof with recursion.
	
	(3) Follows directly from (1) and (2).
\end{proof}

The smallest group that is not totally \semunram is $\SL_2(\FF_3)$ of order $24$. 
The smallest examples of groups that are totally \semunram but not totally \unram are  four groups of order $48$ with GAP group 
IDs $[48,15],[48,16],[48,17]$ and $[48,18]$.
There exist $p$-groups (for every prime $p$) that are not totally \semunram. 
We provide examples in Subsection \ref{sub-p}.

The next example shows that  (normal) subgroups and quotients of a totally \semunram group are not necessarily totally \semunram.
\begin{example}\label{ex-prot}
	The group $G$ with the structure description $C_2 \times((C_8\rtimes C_4)\rtimes C_2)$ and group ID $[128,254]$ is totally \semunram and  the group $H$ with the structure description $(C_8\rtimes C_4)\rtimes C_2)$  and group ID $[64,10]$  is not totally \semunram. 
	However, we have $G\cong C_2\times H$, hence the quotient $G/C_2\cong H$ and the subgroup $1\times H$ of $G$ are not totally \semunram. 
\end{example}

Totally \semunram groups behave well under the direct product.

\begin{proposition}\label{prop-psquo}
	 (1) Let $G$ and $H$ be finite groups and suppose  $G$ is \semunram over $N$, then $G\times H$ is \semunram over $N\times 1$.
	
	(2) If $G$ and $H$ are totally \semunram groups, then $G\times H$   is totally \semunram.
\end{proposition}
\begin{proof}
(1) Any irreducible representation of $N\times 1$ is of the form $\mu\otimes \tau_1$ where $\mu\in \Irr(N)$. 
If $\chi\in \Irr_\mu(G)$ has multiplicity free restriction $\chi|_N$, 
then $\chi \otimes \tau_H \in \Irr_{\mu  \otimes \tau_1}(G \times H)$ has multiplicity free restriction 
$(\chi \otimes \tau_H)|_{N\times 1}=\chi|_N \otimes \tau_1$.

(2) The base case is $A \times H$ where $A$ is abelian. 
By (1), $A\times H$ is \semunram over $A\times 1$ and $(A \times H)/(A\times 1) \cong H$ is totally \semunram. 
We reduce the general case to  the base case recursively using (1).
\end{proof}

We return to Example \ref{ex-prot} and consider it from the point of view of the noncommutative rational invariants.
\begin{example}\label{ex-inv}
	We have $G\cong C_2\times H$ where $G$ is the group with GAP group ID $[128,254]$ that is totally \semunram and $H$ is the group with GAP group ID $[64,10]$ that is not totally \semunram. 
	Let $G$ act on the free skew-field $\frats{x}{n}$ via a complete representation. 
	Then the skew-field of rational invariants  $\frats{x}{n}^G\cong \frats{y}{|G|(n-1)+1}$ is rational by \cite[5.1]{kl20}. 
	However, we can take an intermediate step
	and  first compute the skew-field of $C_2$-invariants $\frats{x}{n}^{C_2}\cong \frats{z}{2n-1}$ that is rational, 
	hence the skew-field
	$\frats{z}{2n-1}^H\cong \frats{x}{n}^G$ is also rational.
	Thus we have an example of a skew-field of rational invariants of a finite group that is rational yet the group is not totally \semunram.
	Furthermore,  $G$ is \semunram over $C_2$, by Proposition \ref{prop-psquo}, 
	thus tracing trough the proofs of Theorem \ref{thm-comppseu} and \cite[5.1]{kl20} we can show that the considered action of $H$ is linear. 
\end{example}

We proceed with a cohomological characterisation of the notion of "\semunram over".

\begin{theorem}\label{thm-mast2}
	Let $N$ be a nontrivial abelian normal subgroup of $G$. 
	Then $G$ is \semunram over $N$ if and only if for all  irreducible characters $\mu\in\Irr(N)$ and 
	for any cocycle $f\in  Z^2(I_G(\mu),N)$ defining the extension from  $I_G(\mu)/N$ to $I_G(\mu)$, 
	the class $[\mu \circ f]\in M(I_G(\mu)/N)$ is trivial.
\end{theorem}
\begin{proof}
	By  (2) of Proposition \ref{prop-lin}, group $G$ is \semunram over $N$ if and only if $\IRR_\mu(I_G(\mu))$ contains a representation of degree $1$. 
		By Lemma \ref{lem-deg}, instead of the degrees of $\IRR_\mu(I_G(\mu))$ we can consider the degrees of irreducible projective representations $\IRR^\pi(I_G(N)/N)$, where $\pi$ is any representative of $[\mu\circ f]$.
	By Lemma \ref{lem-one}, $\IRR^\pi(I_G(N)/N)$  contains a representation of degree $1$ if and only if $[\mu\circ f]$ is trivial.
\end{proof}

We can reformulate  Theorem \ref{thm-mast2} using $T_H$ from the exact sequence \eqref{ex-lhs}.

\begin{corollary}\label{cor-thm2}
A group	$G$ is \semunram over a nontrivial abelian normal subgroup  $N$ if and only for every subgroup $H$ of $G$ containing $N$  the map $T_H\colon  \Lin^H(N)\to M(H/N)$ is trivial. 
	
\end{corollary}
\begin{proof}
	Suppose that $G$ is \unram over $N$. If the subgroup $H$ is of the form $I_G(\mu)$ for some $\mu\in\Irr(N)$, we can directly use Theorem \ref{thm-mast2}. 
	Now take an arbitrary subgroup $H$ of $G$ containing $N$ any $\mu\in\Lin^H(N)$, then a direct computation shows $H\subseteq I_G(\mu)$. 
	It remains to show that $T_H(\mu)$ is trivial, which is true since $T_H(\mu)$ is equal to the restriction of  $T_{I_G(\mu)}(\mu)$ to $H/N$.

The backwards implication follows from applying the assumptions to the inertia subgroups.
\end{proof}

We give two more corollaries to Theorem \ref{thm-mast2}.
\begin{corollary}\label{cor-semdir}
	A semidirect product $A\rtimes G$ of a totally \semunram group $G$ with an abelian group $A$ is totally \semunram.
\end{corollary}
\begin{proof}
		The inertia subgroups of the characters of $A$ are of the form $A\rtimes H$, for some subgroup $H$ of $G$.  
		The $2$-cocycle defining the extension from $H$ to $A\rtimes H$ is trivial, hence  $A\rtimes G$ is \semunram over $A$ by Theorem \ref{thm-mast2}.
		The quotient $(A\rtimes G)/A\cong G$ is totally \semunram, hence $A\rtimes G$ is totally \semunram as well.
\end{proof}
Example \ref{ex-prot} shows that a semidirect product of an abelian group with a totally \semunram group need not be totally \semunram.

\begin{corollary}\label{prop-nuncom2}
	If the commutator subgroup $G'$ is abelian and $G'\neq[G,G,G]$, then $G$ is not \semunram over $G'$. 
\end{corollary}
\begin{proof}
	We consider the exact sequence \eqref{ex-lhs}.
	The invariant characters $\Lin^G(G')\cong G'/[G',G]$ are nontrivial and
	the restriction $\rest\colon \Lin(G) \to \Lin^G(G')$ in \eqref{ex-lhs} is trivial, therefore the map $T_G\colon\Lin^G(G')\to M(G/G')$ is  nontrivial. 
	We now apply Corollary \ref{cor-thm2}.
\end{proof}

\section{Totally unramified groups}\label{sec-unr}

In this section we focus on totally \unram groups.
As mentioned in the previous section
symmetric groups $S_3$ and $S_4$   and all dihedral groups are totally \unram.
The smallest example of a group that is not totally \unram is $\SL_2(\FF_3)$, same as in the totally \semunram case.

 We draw attention to the  theorem that addresses a concept similar to "unramified over".
\begin{theorem}[\cite{myr75}]	Let $H$ be a subgroup of $G$. 
	For every irreducible character $\chi\in \Irr(G)$ the restriction $\chi|_H$ is multiplicity free if and only 
	if the centralizer $C_{\C G}(\C H)$ of the group algebra  $\C H$ is commutative.
\end{theorem}
We get the following corollary.
\begin{corollary}
	Let $N$ be an abelian normal subgroup of $G$. If the centralizer $C_{\C G}(\C H)$ is commutative, then $G$ is \unram over $N$.
\end{corollary}
The converse does not hold, the reason being that in the definition of "\unram over" we allow restrictions to be multiples of the trivial character. 
Any totally \unram group that is not metabelian (such as $S_4$) provides a concrete example where the converse fails.

As in the case  of totally \semunram  groups a (normal) subgroup of a totally \unram group is not necessarily totally \unram. 
\begin{example}
The group $D_{16}\times S_3$ (with GAP group ID $[96,117]$) is totally \unram and contains a normal subgroup with structure description$(C_3\times D_8)\rtimes C_2$ (group ID $[48,15]$), that is not totally \unram.
\end{example}

Totally \semunram groups are closed under direct products but not closed under quotients. In the case of totally \unram groups the things are reversed.
\begin{proposition}\label{prop-quot}	
	Let $\varphi\colon G \to H$ be a surjective group homomorphism.
	\begin{enumerate}[(1)]
		\item If $G$ is \unram over a nontrivial abelian normal subgroup $N$ and $M=\varphi(N)$ is nontrivial, then $H$ is \unram over $M$.
		\item If $G$ is totally \unram, then $H$ is totally \unram.
	\end{enumerate}
	
\end{proposition}
\begin{proof}
	(1)	There are bijections given by inflation  
	$\infl\colon\Irr(H)\to \{\chi \in \Irr(G) \mid \ker\varphi\subseteq\ker\chi\}$
	and 
	$\infl\colon\Irr(M) \to\{\chi \in \Irr(N) \mid \ker\varphi|_N\subseteq\ker\chi\}$. 
	We note that $\infl(\chi|_M)=(\infl\chi)|_N$.

	For any $\chi \in \Irr(H)$ we consider its inflation $\infl\chi \in \Irr(G)$. 
	The restriction $(\infl\chi)|_N$ is multiplicity free or a multiple of the trivial character, whence the same must hold for $\infl(\chi|_M)$ and  $\chi|_M$.

	(2) 	Let $G$ be totally \unram and
	$$
	1= N_0   \subsetneq N_1 \subsetneq \dots\subsetneq N_{n-1} \subsetneq N_n=G
	$$
	a series of normal subgroups such that $N_{j+1}/N_{j}$ is abelian and $G_1/N_{j}$ is \unram over $N_{j+1}/N_{j}$.
	Consider the series
	$$
	1= \varphi(N_0)   \subseteq \varphi(N_1) \subseteq \dots\subseteq \varphi(N_{n-1}) \subseteq \varphi(N_n)=H.
	$$
	We can assume that all the inclusions are strict, otherwise we can remove the redundant elements. 
	By $\bar\varphi\colon G/N_j \to H/\varphi(N_j)$ we denote the induced homomorphism and note that it is surjective. 
	The group $\varphi(N_{j+1})/\varphi(N_j)\cong \bar\varphi(N_{j+1}/N_{j})$  is abelian.
	The group $H/\varphi(N_j)\cong \bar\varphi(G/N_j)$ is \unram over
	$\varphi(N_{j+1})/\varphi(N_{j})\cong \bar\varphi(N_{j+1}/N_{j})$ by  (1).
\end{proof}

An example of a direct product of totally \unram groups that is not totally \unram is $S_4\times S_4$.

The next results narrow down the candidates for abelian normal subgroups $N$ of $G$ over which $G$ is potentially \unram.
\begin{proposition}\label{prop-cent}
	Let $G$ be a finite group and $N$ a nontrivial central subgroup $(1 \subsetneq N \subseteq Z(G))$. 
	Then  $G$ is \unram over $N$ if and only if $G$ is abelian.
\end{proposition}
\begin{proof}
	Because $N$ is central   we have $I_G(\mu)=G$ for every $\mu \in \Irr(N)$. 
By  (1) of Proposition \ref{prop-lin},  we have $\Irr_\mu(G) \subset \Lin(G)$ for any nontrivial $\mu \in \Irr(N)$. 
		Assume that $G$ is not abelian, then there exists $\chi \in \Irr(G)$  with $\chi(1)>1$. We note that $\chi\in\Irr_\tau(G)$.
	Next take any $\theta \in  \Irr_\mu(G)$ for some nontrivial $\mu\in\Irr(N)$. 
	Then $\theta\chi$ is an irreducible  character of $G$, $(\theta\chi)(1)=\chi(1)>1$ and   $(\theta\chi)|_N=\chi(1)\mu$, which is a contradiction.
	
	 The	backwards implication is obvious.
\end{proof}

\begin{corollary}\label{prop-kom2}
	If  a  group $G$ is \unram over $N$, then either $N=[N,G]$ or $[N,G]=G'$.  	In particular $N\subseteq G'$ or $G'\subseteq N$. 
\end{corollary}
\begin{proof}
	If $[N,G] \neq N$, then $G/[N,G]$ is \unram over $N/[N,G]$ by Proposition \ref{prop-quot}. 
	The subgroup $N/[N,G]$ is central in $G/[N,G]$,
	hence $G/[N,G]$ is abelian by Proposition \ref{prop-cent}. 
\end{proof}

To reduce the number of recursive steps in the definition of a totally \unram group we would like to take an abelian normal subgroup as big as possible, however,
there are limitations.
\begin{proposition}\label{prop-ekquo}
	Let $G$ be \unram over $N$ and let $M\supsetneq N$ be an abelian normal subgroup of $G$. 
	Then $G$ is \unram over $M$ if and only if $G/N$ is \unram over $M/N$.
\end{proposition}
\begin{proof}
	Assume $G/N$ is \unram over $M/N$.
	For any $\chi \in \Irr(G)$ we have two options: $\chi|_N$ is a sum of distinct irreducible characters or $\chi|_N$ is a multiple of the trivial character. 
	In the first case $\chi|_M$ must be a sum of distinct irreducible characters of $M$ otherwise an irreducible character would also
	appear  multiple times in  $(\chi|_N)=(\chi|_M)|_N$.
	In the second case $\chi$ is induced from a character $\bar\chi\in\Irr(G/N)$. 
	We have $\chi|_M=\Ind(\bar\chi|_{M/N})$.
	Since 
	$\bar\chi|_{M/N}$ is  a multiple of the  trivial character or multiplicity free, so is $\chi|_M$.
	
	Conversely, if $G$ is \unram over $M$, then  by Proposition \ref{prop-quot}, $G/N$ is \unram over $M/N$.
\end{proof}

\begin{corollary}\label{cor-ext}
	Let $G$ be \unram over $N$ and assume $Z(G)\not\subset N$. 
	Then $G$ is \unram over $Z(G) N$  if and only if $G'\subseteq N$.
\end{corollary}
\begin{proof}
	
	By Propositions \ref{prop-cent},  $G$ is  \unram over $Z(G)N$ if and only if $G/N$ is \unram over $Z(G)N/N$. 
	We notice that $Z(G)N/N$ is nontrivial and central in $G/N$, therefore, by Proposition \ref{prop-ekquo},
	$G$ is \unram over $Z(G)N$ if and only if $G/N$ is abelian.
\end{proof}

There indeed exist  examples of totally \unram groups that are not \unram over any abelian normal subgroup that contains the center. 
One such group is $C_3\rtimes((C_{10}\times C_2)\rtimes C_2)$  with GAP group ID $[120,11]$.

We proceed with a cohomological characterisation of the notion of "\unram over".

\begin{theorem}\label{thm-mast}

Let $N$ be a nontrivial abelian normal subgroup of $G$.
	Then $G$ is \unram over $N$ if and only if for all nontrivial $\mu \in \Irr(N)$ the following conditions are satisfied:
	\begin{enumerate}
		\item $I_G(\mu)/N$ is abelian;
		\item for any  $2$-cocycle $f\in Z^2(I_G(\mu), N)$ defining the extension from  $I_G(\mu)/N$ to $I_G(\mu)$, the class of the factor set  $[\mu \circ f]\in M(I_G(\mu)/N)$ is trivial.
	\end{enumerate}
\end{theorem}
\begin{proof}
	
	By (1) of Proposition \ref{prop-lin}, $G$ is \unram over $N$ if and only if for every nontrivial character $\mu \in \Irr(N)$, every representation in $\IRR_\mu(I_G(\mu))$ is of degree one. 
	By Lemma \ref{lem-deg}, instead of the degrees of $\IRR_\mu(I_G(\mu))$ we can consider the degrees of irreducible projective representations $\IRR^\pi(I_G(N)/N)$, where $\pi$ is any representative of $[\mu\circ f]$.

	If $[\mu \circ f]$ is trivial, then we   consider the degrees of linear representations $\IRR(I_G(\mu)/N)$. 
	We note that $I_G(\mu)/N$ is abelian if and only if all its irreducible linear representations are of degree one. 	
	If $[\mu \circ f]$ is not trivial, then by Lemma \ref{lem-one}, none of  the irreducible $(\mu \circ f)$-representations of $I_G(\mu)/N$ is of degree one.
 Hence $\IRR_\mu(I_G(\mu))$ contain only representations of degree $1$ if and only if $I_G(\mu)/N$ is abelian and  $[\mu\circ f]$ trivial.
\end{proof}

Next we combine Theorem \ref{thm-mast} with the exact sequence \eqref{ex-lhs}.
\begin{corollary}
	A group $G$ is \unram over $N$ if and only if for every  subgroup $H$ of $G$ containing $N$ such that $\Lin^H(N)$  is nontrivial, the map $T_H\colon \Lin^H(N)\to M(H/N)$  is trivial and $H/N$ is abelian.
\end{corollary}
\begin{proof}
We begin with the forward implication.	If $H$ is of the form $I_G(\mu)$ for some $\mu\in\Irr(N)$, we can directly use Theorem \ref{thm-mast}. 
Next take an arbitrary subgroup $H$ containing $N$ and a nontrivial character $\mu\in\Lin^H(N)$, then $H\subseteq I_G(\mu)$, therefore $H/N$ is abelian. 
Also $T_H(\mu)$ is equal to the restriction of  $T_{I_G(\mu)}(\mu)$ to $H/N$, hence it is trivial.
	
	To prove the  backwards implication we just apply the assumptions to the inertia subgroups.
\end{proof}
Using Theorem \ref{thm-mast} we provide two classes of totally \unram groups.

\begin{corollary}\label{prop-cikli}
	(1) A group $G$ is \unram over an abelian normal subgroup $N$ if $I_G(\mu)/N$ is cyclic for all $\mu \in \Irr(N)$.
	In particular, if $G/N$ is cyclic, then $G$ is \unram over $N$.

(2) Metacyclic groups are totally \unram. 
\end{corollary}
\begin{proof}
Since the Schur multiplier of a cyclic group is trivial  we can apply Theorem \ref{thm-mast}.
\end{proof}

\begin{corollary}\label{prop-semd}
	Semidirect products of abelian groups are totally \unram.
\end{corollary}
\begin{proof}
	Let $G=A\rtimes B$ be a semidirect product of abelian groups. We show that $G$ is \unram over $A$. The inertia subgroups of the characters of $A$ are of the form
	  $A\rtimes C$, for some subgroup $C$ of $B$. 
	 Then the $2$-cocycle $f$ defining the extension from $A$ to $A\rtimes C$ is associated to the trivial one  and hence $\la \circ f$ is associated to the trivial $2$-cocycle.
	 We now apply Theorem \ref{thm-mast}.
\end{proof}

\subsection{Isoclinism} The notion of isoclinism was introduced by Hall \cite{hall1940}.
For any group $G$  we have the induced commutator map 
$ G/Z(G)\times G/Z(G) \to G'$
defined by 
$$
(g_1Z(G),g_2Z(G)) \mapsto [g_1,g_2].
$$
Groups $G$ and $H$ are \emph{isoclinic} if there exist isomorphisms
$\varphi \colon G/Z(G) \to H/Z(H)$ 
and 
$\psi\colon G' \to H'$
that commute with the commutator map, 
that is, if $\varphi(g_iZ(G))=h_iZ(H)$ for $i=1,2$, then 
$\psi[g_1,g_2] =[h_1,h_2]$.
An equivalence class with respect to isoclinism is called  an \emph{ (isoclinism) family}, a group of the smallest order in a family  is called  \emph{a stem group}.

Groups from the same family share some properties concerning their representations and we will make use of this in our study of totally \unram groups.

\begin{lemma}\label{lem-iso1}
	Let $G$ and $H$ be finite groups and let $\varphi\colon G/Z(G) \to H/Z(H)$ and $\psi\colon G' \to H'$ satisfy the isoclinism conditions and let $N$ be a nontrivial normal subgroup of $G$ containing $Z(G)$ and let $M$ be the biggest normal subgroup of $H$ satisfying $\varphi(N/Z(G))=\varphi(M/Z(H))$.
	\begin{enumerate}[(1)]
		\item If $N$ is abelian, then $M$ is abelian.
		\item $G/N$ is isomorphic to $H/M$.
		\item Let $\rho_H$ be an irreducible linear representation of $H$, $P$ a descent of $\rho$ to $H/Z(H)$ and $P_G$ any lift of $\varphi^*\rho$ to a linear representation of $G$. 
		The restriction $\rho_H|_M$ is multiplicity free (trivial) if and only if $\rho_G|_N$ is multiplicity free (trivial).
		\item If $N$ is abelian,  then $G$ is \unram over $N$ if and only if $H$ is \unram over $M$.
	\end{enumerate}  
\end{lemma}
\begin{proof}
	(1)	For any $h_1,h_2 \in M$  take $g_1,g_2\in N$  such that $\varphi(g_i Z(G))= h_iZ(H)$,
	then $[h_1,h_2]=\psi[g_1,g_2]=\psi(1)=1$. 
	
	(2) By Noether's isomorphism theorems we have $$G/N\cong \frac{G/Z(G)}{N/Z(G)} \cong \frac{H/Z(H)}{M/Z(H)}\cong H/M.$$
	
	(3) It is enough to prove the statement in one direction. The other follows from the symmetry of isoclinism.
	Let $\rho_H$ be an irreducible linear representation of $H$ of degree $n$ and let  $\rho_H|_M=\bigoplus_{i=1}^{n} \sigma_H^i$ be the decomposition into irreducible summands. 
	Suppose the restriction is multiplicity free or trivial.
	Let $P$ be the projective representation obtained by descent of $\rho$ to $H/Z(H)$.
	Then $P|_{M/Z(H)}$ is a direct sum of descents of $\sigma_H^i$, $i=1,\dots,n$.
	Let the linear representation $\rho_G$ of $G$ be a lift of $\varphi^*P$.
Then  $\varphi^*(P|_{M/Z(H)})$  is also descent of $\rho_G|_N$, 
hence  $\rho_G|_N$ is multiplicity free (trivial) as well.
	
	(4)	
	By \cite[Ch.3, Cor.2.5]{beyl2006group}, an  irreducible projective representation of $G/Z(G)\cong H/Z(G)$ lifts to a linear representation of $G$ 
	if and only if it lifts to a linear representation of $H$, therefore, we can use  (3) on each irreducible representation of $H$.
\end{proof}

We move to a case where $G$ is \unram over a normal abelian subgroup $N$ with $N\subset G'$.

\begin{lemma}\label{lem-iso2}
	Let $G$ and $H$ be finite groups and let $\varphi\colon G/Z(G) \to H/Z(H)$ and $\psi\colon G' \to H'$ satisfy the isoclinism conditions.  
	
	\begin{enumerate}[(1)]
		\item For any $g\in G$ and $h\in H$, such that $\varphi(gZ(G))=hZ(H)$, and any $n \in G'$ we have $\psi(gng\inv)=h\psi(n)h\inv$.
		In particular, for any normal subgroup $N$ of $G$ contained in $G'$ the image $\psi(N)$ is a normal  subgroup of $H$. 		
		\item 
		For any normal subgroup $N$ of $G$ contained in $G'$ the quotient  $G/N$ is isoclinic to $H/\psi(N)$.	
		
		\item For any projective representation $P$ of $H/Z(H)$, its lift to a linear representation $\rho_H$ of group  the $H$ and any lift  of a  representation $\varphi^*P$ to a linear representation $\rho_G$ of group the $G$, we have $\rho_G|_N=\psi^*(\rho_H|_M)$.
		\item If $G$ in \unram over a nontrivial normal abelian subgroup $N$ contained in $G'$ then $H$ is unramified over $\psi(N)$.
	\end{enumerate}
\end{lemma}
\begin{proof}
	(1) It is enough to consider commutators. Let  
	$n=[g_1,g_2]$ 
	and $m=\psi(n)$.
	Take $h_1,h_2\in H$  such that 
	$\varphi(g_i Z(G))=h_iZ(H)$,
	then we have
	$m=[h_1,h_2]$.
	Take any $h\in H$ and let $g \in G$ be such that $\varphi(gZ(G))=hZ(H)$.
	We get
	$$
	hm h\inv  =[hh_1h\inv,hh_2h\inv]
	\quad \text{and} \quad
	gn g\inv  =[gg_1g\inv,gg_2g\inv],
	$$
	which shows  that $\psi(gng\inv)=hmh\inv$.
	
	(2)  Denote $M=\psi(N)$.
	For $(gN)Z(G/N)\in (G/N)/Z(G/N)$ we  define 
	$\bar\varphi ((gN)Z(G/N)) = (hM)Z(H/M),$
	where $h\in H$ is any element satisfying $\varphi(g Z(G))=hZ(H)$. 
	This yields a well defined isomorphism. 
	Then $\bar\varphi\colon (G/N)/Z(G/N) \to (H/M)/Z(H/M)$ and $\bar\psi \colon G'/N \to H'/M$, where $\bar\psi(xN)=\psi(x)M$,  satisfy the isoclinism conditions.
	
	(3) Pick transversals $\{h_\alpha\mid \alpha \in H/Z(H)\}$ and 
	$\{g_\alpha\mid \alpha \in G/Z(G)\}$. 
	Representations $\rho_H$ and $\rho_G$ are given by 
	$\rho_H(zh_\al)=\la(z)P(\al)$ and $\rho_G(wG_\al)=\mu(w)P(\varphi(\al))$,
	for some linear characters $\la$ and $\mu$ of $Z(H)$ and $Z(H)$, respectively

	It is enough to show the desired equality $\rho_G(n)=\psi^*(\rho_H|_M)(n)$ for a
	commutator
	$
	n=[g_\alpha,g_\beta]
	$
	where  $g_\alpha$ and $g_\beta$ are  from the transversal. 
	We show the identity by expanding
	$$
	\rho_G(n)=[\rho_G(g_\al),\rho_G(g_\beta)]=[P(\varphi(\alpha)),P(\varphi(\beta))]
	$$
	and
	$$
	\rho_H(\psi(n))=\rho_H([h_{\varphi(\alpha)}, h_{\varphi(\beta)}])=[P(\varphi(\alpha)),P(\varphi(\beta))].
	$$

	(4) 
	Any irreducible linear representation of $G$ is a lift of an irreducible projective representation of $G/Z(G)$. 
	By \cite[Ch.3, Cor.2.5]{beyl2006group}, an  irreducible projective representation of $G/Z(G)\cong H/Z(G)$ lifts to a linear representation of $G$ 
	if and only if it lifts to a linear representation of $H$. 
	The restriction of lifts to $G'$ and $H'$, respectively,  yield equivalent representations of, by  (3), and therefore, restrictions to $N\cong\psi(N)$ are equivalent as well. 
	Hence, the restriction of every irreducible linear representation of $G$ to $N$ is multiplicity free or trivial 
	if and only if 
	 restriction of every irreducible linear representation of $H$ to $\psi(N)$ is multiplicity free or trivial. 
\end{proof}

Finally we combine the two cases.

\begin{theorem}\label{prop-piso}
	If an isoclinic family contains a totally \unram group, then every finite group in this family is totally  \unram.
\end{theorem}
\begin{proof} 
	All finite groups in the family of abelian groups are totally \unram.
	Consider  finite representatives  $G$ and $H$ of a nonabelian isoclinic family.
	If $G$ is totally \unram then there exists a normal abelian group $N$ such that $G$ is \unram over $N$, $G/N$ is totally unramified and  $G'\subseteq N$ or $N \subseteq G'$ by Corollary \ref{prop-kom2}. 
	In the case  $G'\subseteq N$ we can always assume $Z(G)\subseteq N$,  otherwise use Corollary \ref{cor-ext} to replace $N$ with $Z(G)N$. 
	Then  we use Lemma \ref{lem-iso1} to  find a normal abelian subgroup $M$ of $H$ such that $H$ is \unram over $M$ and $H/M\cong G/N$ is abelian, hence, totally \unram.
	In the case  $N \subseteq G'$ we use Lemma \ref{lem-iso2} to find a normal abelian subgroup $M$ of $H$ such that $H$ is \unram over $M$ and $H/M$ is isoclinic to $G/N$ and finish the proof by recursion.
\end{proof}

Although a direct product of totally \unram group need not be totally \unram we have the following weaker result.
\begin{corollary}\label{prop-dir}
	A direct product of a totally \unram group and an abelian group is totally \unram.
\end{corollary}
\begin{proof}
A group $G$ is isoclinic to $G\times A$ for any abelian group $A$, hence, the result follows from Theorem \ref{prop-piso}.
\end{proof}

\section{Nilpotent totally \unram groups}\label{sec-nilunr}

In this section we prove some stronger results for nilpotent totally \unram groups.
\begin{corollary}\label{prop-nuncom}
	If $G$ is metabelian nilpotent, then it is not \unram over $G'$. 
\end{corollary}
\begin{proof}
	Directly from  Corollary \ref{prop-nuncom2}.
\end{proof}
This yields a new proof of a well-known result.
\begin{corollary}\label{cor-necic}
	If  $G$ is a metabelian nilpotent group that is not abelian, then $G/G'$ is not cyclic.
\end{corollary}
\begin{proof}
	By Corollary \ref{prop-nuncom}, $G$ is not \unram over $G'$, hence
	 the quotient $G/G'$ is not  cyclic by Proposition \ref{prop-cikli}.
\end{proof}

We can improve Corollary \ref{prop-kom2} for nilpotent groups.
\begin{corollary}\label{prop-kom3}
	If  a  nilpotent group  $G$ is \unram over $N$, then $[N,G]=G'$.  
\end{corollary}
\begin{proof}
	Follows directly from Corollary \ref{prop-kom2}, as a  nilpotent group $G$  does not contain any nontrivial abelian subgroup $N$ with  property $[N,G]=N$. 
\end{proof}

Being nilpotent totally \unram is a very restrictive property. 

\begin{theorem}\label{prop-com}
	If a nilpotent group $G$ is \unram over an  abelian normal subgroup $N$, then $G'\subsetneq N$. 
	In particular every totally \unram nilpotent group is metabelian.
\end{theorem}
\begin{proof}
	 By Corollary \ref{prop-kom3}, we get $G' \subseteq N$ and 
	the equality is excluded since $[G',G] \subsetneq G'$.
\end{proof}

Contrary to the general case   nilpotent totally \unram groups are closed under the direct product.
\begin{proposition}\label{prop-prod}
	If $G_1$ and $G_2$ are totally \unram  nilpotent groups, then $G_1\times G_2$ is totally \unram.
\end{proposition}
\begin{proof}
	Let $G_i$ be \unram over $N_i$, we show that $G_1\times G_2$ is \unram over $N_1\times N_2$, then the quotient is abelian and we are done.
	
	Every irreducible character of $G_1 \times G_2$ is of the form $\chi_1 \otimes \chi_2$ where $\chi_i$ is an irreducible character of $G_i$. 
	We get $(\chi_1 \otimes \chi_2)|_{N_1 \times N_2}= \chi_1|_{N_1} \otimes \chi_2|_{N_2}$
	and each $\chi_i|_{N_i}$ is multiplicity free or a multiple of the trivial character.
	We consider two cases; either both of $\chi_i|_{N_i}$ are multiplicity free or at least one of them is a multiple of  the trivial character.
	
	Let 
	$\chi_1|_{N_1} = \sum_{i=1}^m \mu_i$
	 and 
	 $\chi_2|_{N_2} = \sum_{j=1}^n \theta_j$,
	  where $\mu_i$ and $\theta_j$  are pairwise distinct characters of $N_1$ and $N_2$ respectively.
	 Then 
	$\chi_1|_{N_1}\otimes \chi_2|_{N_2} = \sum_{i,j} \mu_i \otimes \theta_j$ and 
	$\la_i \otimes \chi_j$ are pairwise distinct.
	
	Let $\chi_1$ be a multiple of the trivial character. 
	By Proposition \ref{prop-com}, $N_1$ contains $G_1'$ and therefore $\ker(\chi_1)$ contains $G_1'$. This forces   $\chi_1$ to be linear. 
	The restriction $\chi_1|_{N_1} \otimes \chi_2|_{N_2}=\tau \otimes \chi_2|_{N_2}$ is clearly multiplicity free  or a multiple of the trivial character. 
	The case where $\chi_2|_{N_2}$ is a multiple of the trivial character is  symmetric.
\end{proof}

If $G$ is \unram over $N$ and $M$ is an abelian normal subgroup containing $N$,  $G$ is not necessarily \unram over $M$. 
This obstacle does not apply to nilpotent groups.
\begin{proposition}\label{prop-pmax}
	If a nilpotent group $G$ is \unram  over a nontrivial abelian normal subgroup $N$, then $G$ is \unram over any abelian normal subgroup containing $N$.
\end{proposition}
\begin{proof}
	Let $M$ be an abelian normal subgroup of $G$ containing $N$. 
	By Theorem \ref{prop-com}, $G/N$ is abelian. An abelian group is \unram over any nontrivial subgroup. 
	By Proposition \ref{prop-ekquo}, $G$ is \unram over $M$.
\end{proof}

\subsection{Totally \unram \texorpdfstring{$p$}{TEXT}-groups}\label{sub-p}
Every nilpotent group is a direct product of $p$-groups. 
Thus  Propositions \ref{prop-quot} and \ref{prop-prod} show that a nilpotent group is totally \unram if and only if its Sylow $p$-subgroups are totally \unram. It is therefore of interest to understand totally \unram $p$-groups.
We  classify totally \unram $p$-group of rank at most $5$
starting with rank at most $4$.

\begin{proposition}\label{prop-4}
	Every $p$-group of rank at most $4$ is totally \unram.
\end{proposition}
\begin{proof}
	Groups of order $p$ and  $p^2$ are abelian and every group of order $p^3$ is metacyclic, hence totally \unram by Corollary \ref{prop-cikli}.
	Every group $G$ of order $p^4$ has an abelian normal subgroup  $N$ of order $p^3$ (\cite{bur13}). 
	Then $G/N$ is cyclic and by Proposition \ref{prop-cikli}, $G$ is totally \unram.	
\end{proof}

We continue with groups of order $p^5$.
We separate the cases $p=2$ and $p\geq3$.
For $p=2$ we use GAP for direct computation. All groups of order $2^5$ are totally \unram. 
There are groups of order $2^6$ that are not totally \unram but are totally \semunram and also groups that are not totally \semunram.

To classify the totally \unram $p$-groups it is enough, by Proposition \ref{prop-piso}, to classify isoclinism families that contain a totally \unram group. 
For $p\geq 3$ we use the classification of $p$-groups by \cite{Jam80}, where we also refer to for the explanation of the classification and notation. 
We just mention that isoclinism families are denoted by $\Phi_s$, $s=1,2,\dots$ and 
if the word  $[g_i,g_j]$ where $g_i$ and $g_j$ are generators does not appear among relations of the group presentation, the relation $[g_i,g_j]=1$ should be assumed.

There are 10  families of $p$-groups ($p\geq 3$) that contain a stem group of rank at most $5$; $\Phi_1,\Phi_2,\dots,\Phi_{10}$.
The family $\Phi_1$ contains all abelian groups.
The families $\Phi_2$ and $\Phi_3$ have a stem group of rank $3$ and $4$  respectively which are totally \unram  by Proposition \ref{prop-4}. 
The remaining families with a stem group of rank $5$ are $\Phi_4,\Phi_5,\dots, \Phi_{10}$. 
We deal with each family separately. 

By Propositions \ref{prop-pmax} and Theorem \ref{prop-com}, to check whether a $p$-group $G$ is totally \unram  we only have to consider maximal abelian normal subgroups of $G$ that strictly contain $G'$. 
We use Corollaries \ref{prop-cikli} and \ref{prop-semd} to give positive answers.

\begin{enumerate}[(i)]
	\item[($\Phi_4$)]
	We consider
	$$
	\Phi_4(1^4)=\langle
	\alpha,\alpha_1,\alpha_2, \beta_1, \beta_2 
	\mid
	[\alpha_i,\alpha]=\beta_i, \alpha^p=\alpha_i^p=\beta_i^p=1 \text{ for } i=1,2
	\rangle.
	$$
	The subgroup 
	$\langle \alpha, \beta_1,\beta_2\rangle $
	is abelian normal and has 
	the abelian subgroup 
	$\langle \alpha_1, \alpha_2 \rangle$
	as a complement, 
	i.e.,  the group $\Phi_4(1^4)$ is a semidirect product of abelian groups and therefore totally  \unram.

	\item[($\Phi_5$)]
	We consider
	$$
	\Phi_5(1^5)=
	\langle
	\alpha_1,\alpha_2,\alpha_3,\alpha_4, \beta 
	\mid
	[\alpha_1, \alpha_2]=[\alpha_2, \alpha_3]=\beta, \alpha_i^p=\beta^p=1 
	\text{ for } i=1,2,3,4
	\rangle.
	$$
	The subgroup 
	$\langle \alpha_1, \alpha_3,\beta\rangle $
	is abelian normal and has 
	the abelian subgroup 
	$\langle \alpha_2, \alpha_4 \rangle$
	as a complement, 
	i.e., the group $\Phi_5(1^5)$ is a semidirect product of abelian groups and therefore totally \unram.
	
	\item[($\Phi_6$)] 
	We consider 
	$$
	\Phi_6(1^5)=
	\langle
	\alpha_1, \alpha_2, \beta , \beta_1, \beta_2
	\mid
	[\alpha_1,\alpha_2]=\beta, [\beta, \alpha_i]=\beta_i, \alpha_i^p=\beta^p=\beta_i^p=1 
	\text{ for } i=1,2
	\rangle.
	$$
	The commutator subgroup 
	$\langle \beta, \beta_1, \beta_2\rangle$
	is a maximal abelian normal subgroup,
	therefore  the group $\Phi_6(1^5)$  is not totally \unram.
	
	\item[($\Phi_7$)]
	We consider 
	\begin{equation*}
	\begin{aligned}
	\Phi_7(1^5)=
	\langle
	\alpha,\alpha_1,\alpha_2,\alpha_3,\beta 
	\mid
	[\alpha_i,\alpha]=\alpha_{i+1} ,[\alpha_1,\beta]=\alpha_3, 
	\alpha^p=\alpha_1^{(p)}=\alpha_{i+1}^p=\beta^p=1
	\text{ for } i=1,2
	\rangle.
		\end{aligned}
		\end{equation*}
	The subgroup 
	$\langle \alpha_1, \alpha_2,\alpha_3\rangle $
	is abelian normal and has 
	the abelian subgroup 
	$\langle \alpha, \beta \rangle$
	as a complement, 
	i.e.,  the group $\Phi_7(1^5)$ is a semidirect product of abelian groups and therefore totally \unram.

	\item[($\Phi_8$)]
	We consider 
	$$
	\Phi_8(32)=
	\langle
	\alpha_1, \alpha_2,\beta
	\mid
	[\alpha_1,\alpha_2]=\beta=\alpha_1^p, \beta^{p^2}=\alpha_2^{p^2}=1
	\rangle.
	$$
	The subgroup 
	$\langle \alpha_1, \beta \rangle$
	is abelian normal and has a cyclic quotient 
	$\langle \alpha_2 \rangle$. Therefore
	the group $\Phi_9(1^5)$ is  totally \unram.

	\item[($\Phi_9$)]
	We consider 
	$$
	\Phi_9(1^5)=
	\langle
	\alpha, \alpha_1, \alpha_2, \alpha_3,\alpha_4
	\mid
	[\alpha_i , \alpha]=\alpha_{i+1}, \alpha^p=\alpha_1^{(p)}=\alpha_{i+1}^{(p)}=1
	\text{ for } i=1,2,3
	\rangle.
	$$
	The subgroup 
	$\langle \alpha_1, \alpha_2, \alpha_3,\alpha_4 \rangle$ 
	is abelian normal and has a cyclic quotient 
	$\langle \alpha \rangle$. Therefore
	  the group $\Phi_9(1^5)$ is  totally \unram.

	\item[($\Phi_{10}$)]
	We consider 
\begin{equation*}
	\begin{aligned}
	\Phi_{10}(1^5)=
	\langle
	\alpha, \alpha_1, \alpha_2, \alpha_3,\alpha_4
	\mid
	[\alpha_i , \alpha]=\alpha_{i+1}, [\alpha_1 , \alpha_2]=\alpha_4,
	\alpha^p=\alpha_1^{(p)}=\alpha_{i+1}^{(p)}=1
	\text{ for } i=1,2,3
	\rangle.
\end{aligned}
\end{equation*}
	The commutator subgroup 
	$\langle \alpha_2, \alpha_3, \alpha_4 \rangle$
	is a maximal abelian normal subgroup,
	therefore  the group $\Phi_6(1^5)$  is not totally \unram.
	\end{enumerate}

\begin{remark}
	There exist totally \unram groups (of order at least $p^4$)  that are not semidirect products of abelian groups. 
	Additionally, $\Phi_7(2111)$ is totally \unram but is not a semidirect product of abelian groups and does not contain any abelian normal subgroup with a cyclic quotient.
\end{remark}

The classification of totally \unram $p$-groups of rank at most $5$ also classifies  totally \semunram $p$-groups of rank at most $5$.
\begin{proposition}\label{prop-=5}
	A $p$-group of rank  at most $5$ is totally \unram if and only if it is totally \semunram.
\end{proposition}
\begin{proof}
	We only have to prove the backwards implication for $p\geq3$. 
	If $G$ is  $p$-group of rank at most $5$ and not totally \unram, then it is from the isoclinism family $\Phi_{6}$ or $\Phi_{10}$. 
	We note that these groups have  $|G'|=p^3$ and hence  $|\Lin(G)|=|G/G'|=p^2$.
	
	Assume $G$  is \semunram over $N$ and not totally \unram.
	Then $G/N$ is not abelian, otherwise $G$ would be \unram over $N$, 
thus totally \unram. However, $G/N$  is metabelian, therefore 
	$
	|\Lin(G/N)|= |(G/N)/(G/N)'|\geq p^2
	$
	by Corollary \ref{cor-necic}. 
	Then the $\inf\colon \Lin(G/N) \to \Lin(G)$ is surjective. 
	 Restriction and $T_G$  from the exact sequence  \eqref{ex-lhs} are trivial maps,   hence $\Lin^G(N)\cong N/[N,G]$ must be trivial which is a contradiction since $N \neq [N,G]$. 
\end{proof}

We give reasoning for the family $\Phi_{11}$ with a stem group of rank $6$ which provides  examples of groups of nilpotency class 2 that are not totally \unram and an example of a  totally \semunram $p$-group that is not totally \unram.

\begin{enumerate}[(i)]

	\item[($\Phi_{11}$)]
	
	We consider
\begin{equation*}
	\begin{aligned}
	\Phi_{11}(1^6)=
	\langle
	\alpha_1, \beta_1,\alpha_2, \beta_2,\alpha_3, \beta_3
	\mid
	[\alpha_1 , \alpha_2]=\beta_{3}, [\alpha_2 , \alpha_3]=&\beta_1, [\alpha_3,\alpha_1]=\beta_3,\\
	 &\alpha_i^p=\beta_i^{p}=1
	\text{ for } i=1,2,3
	\rangle.
\end{aligned}
\end{equation*}
	The commutator subgroup $\langle  \beta_1,\beta_2, \beta_3  \rangle$ is also the center. 
	There are three maximal abelian normal subgroups $N_i=\langle  \alpha_i, \beta_1,\beta_2, \beta_3  \rangle$ for  $i=1,2,3$. 
	We have $\beta_i \not\in [N_i, G]$, therefore $G$ is not \unram over $N_i$ by Corollary \ref{prop-kom3}, hence not totally \unram.

\end{enumerate}
Furthermore, the decomposition  $\Phi_{11}(1^6) =\langle \al_1,\beta_2,\beta_3\rangle\rtimes \langle \beta_1,\al_2,\al_3 \rangle $  shows that $\Phi_{11}(1^6)$ is totally \semunram by Corollary \ref{cor-semdir}.
Hence, we cannot extend the Proposition \ref{prop-=5} to  groups of order $p^6$.


\subsection*{Acknowledgment}
The author thanks Igor Klep, Urban Jezernik and Primož Moravec for fruitful discussions

\bibliography{references}

\newcommand{\etalchar}[1]{$^{#1}$}
\begin{thebibliography}{CHKK10}

\bibitem[Ami66]{ami66}
SA~Amitsur.
\newblock Rational identities and applications to algebra and geometry.
\newblock {\em Journal of Algebra}, 3(3):304--359, 1966.

\bibitem[Ber90]{ber90}
George~M Bergman.
\newblock Ordering coproducts of groups and semigroups.
\newblock {\em Journal of Algebra}, 133(2):313--339, 1990.

\bibitem[BR11]{ber11}
Jean Berstel and Christophe Reutenauer.
\newblock {\em Noncommutative rational series with applications}, volume 137.
\newblock Cambridge University Press, 2011.

\bibitem[BT06]{beyl2006group}
F.~Rudolf Beyl and J{\"u}rgen Tappe.
\newblock {\em Group extensions, representations, and the Schur multiplicator},
  volume 958.
\newblock Springer, 2006.

\bibitem[Bur13]{bur13}
William Burnside.
\newblock On some properties of groups whose orders are powers of primes.
\newblock {\em Proceedings of the London Mathematical Society}, 2(1):225--245,
  1913.

\bibitem[BZ97]{BE98}
Yakov Berkovich and Emmanuel Zhmud'.
\newblock {\em Characters of finite groups. Part 1}.
\newblock Number pt. 1 in Translations of mathematical monographs. American
  Mathematical Soc., 1997.

\bibitem[CHKK10]{chu10}
Huah Chu, Shou-Jen Hu, Ming-chang Kang, and Boris~E Kunyavskii.
\newblock Noether’s problem and the unramified brauer group for groups of
  order 64.
\newblock {\em International Mathematics Research Notices},
  2010(12):2329--2366, 2010.

\bibitem[Coh95]{Coh95}
Paul~Moritz Cohn.
\newblock {\em Skew Fields: Theory of General Division Rings}.
\newblock Encyclopedia of Mathematics and its Applications. Cambridge
  University Press, Cambridge, 1995.

\bibitem[Coh06]{Coh06}
Paul~Moritz Cohn.
\newblock {\em Free ideal rings and localization in general rings}, volume~3 of
  {\em New Mathematical Monographs}.
\newblock Cambridge University Press, Cambridge, 2006.

\bibitem[Con94]{con94}
Alain Connes.
\newblock {\em Noncommutative geometry}.
\newblock Academic Press, Inc., San Diego, CA, 1994.

\bibitem[CR81]{CR81}
Charles~W. Curtis and Irving Reiner.
\newblock {\em Methods of representation theory}, volume~2.
\newblock Wiley-Interscience, 1981.

\bibitem[CSST20]{cec20gelf}
Tullio Ceccherini-Silberstein, Fabio Scarabotti, and Filippo Tolli.
\newblock {\em Gelfand Triples and Their Hecke Algebras: Harmonic Analysis for
  Multiplicity-Free Induced Representations of Finite Groups}, volume 2267.
\newblock Springer Nature, 2020.

\bibitem[CTS07]{C-TS07}
Jean-Louis Colliot-Th\'el\`ene and Jean-Jacques Sansuc.
\newblock The rationality problem for fields of invariants under linear
  algebraic groups (with special regards to the {B}rauer group).
\newblock In {\em Algebraic groups and homogeneous spaces}, volume~19 of {\em
  Tata Inst. Fund. Res. Stud. Math.}, pages 113--186. Tata Inst. Fund. Res.,
  Mumbai, 2007.

\bibitem[DR97]{duc97}
G{\'e}rard Duchamp and Christophe Reutenauer.
\newblock Un critere de rationalit{\'e} provenant de la g{\'e}om{\'e}trie non
  commutative.
\newblock {\em Inventiones Mathematicae}, 128(3):613--622, 1997.

\bibitem[EM73]{end73}
Shizuo End{\^o} and Takehiko Miyata.
\newblock Invariants of finite abelian groups.
\newblock {\em Journal of the Mathematical Society of Japan}, 25(1):7--26,
  1973.

\bibitem[GAP21]{GAP4}
The GAP~Group.
\newblock {\em {GAP -- Groups, Algorithms, and Programming, Version 4.11.1}},
  2021.

\bibitem[GKL{\etalchar{+}}95]{gel95}
Israel~M Gelfand, Daniel Krob, Alain Lascoux, Bernard Leclerc, Vladimir~S
  Retakh, and Jean-Yves Thibon.
\newblock Noncommutative symmetric functions.
\newblock {\em Advances in Mathematics}, 112:218--348, 1995.

\bibitem[Hal40]{hall1940}
Philip Hall.
\newblock The classification of prime-power groups.
\newblock {\em Journal f{\"u}r die reine und angewandte Mathematik},
  182:130--141, 1940.

\bibitem[HKK13]{hos13}
Akinari Hoshi, Ming-Chang Kang, and Boris~E Kunyavskii.
\newblock Noether's problem and unramified brauer groups.
\newblock {\em Asian Journal of Mathematics}, 17(4):689--714, 2013.

\bibitem[Isa94]{Isa94}
I.~Martin Isaacs.
\newblock {\em Character theory of finite groups}, volume~69.
\newblock Courier Corporation, 1994.

\bibitem[Jam80]{Jam80}
Rodney James.
\newblock The groups of order $p^{6}$ ($p$ and odd prime).
\newblock {\em Mathematics of computation}, pages 613--637, 1980.

\bibitem[Kar87]{Kar87}
Gregory Karpilovsky.
\newblock {\em The schur multiplier}.
\newblock Oxford University Press, Inc., 1987.

\bibitem[KPPV20]{kl20}
Igor Klep, James~Eldred Pascoe, Gregor Podlogar, and Jurij Vol{\v{c}}i{\v{c}}.
\newblock Noncommutative rational functions invariant under the action of a
  finite solvable group.
\newblock {\em Journal of Mathematical Analysis and Applications},
  490(2):124341, 2020.

\bibitem[KVV12]{kal12}
Dmitry~S Kaliuzhnyi-Verbovetskyi and Victor Vinnikov.
\newblock Noncommutative rational functions, their difference-differential
  calculus and realizations.
\newblock {\em Multidimensional Systems and Signal Processing}, 23(1):49--77,
  2012.

\bibitem[Len80]{len80}
HW~Lenstra.
\newblock Rational functions invariant under a cyclic group.
\newblock In {\em Proceedings of the Queen's number theory conference, 1979,
  Queen's papers in pure and applied mathematics}, volume~54, pages 91--99,
  1980.

\bibitem[Lew74]{lew74}
Jacques Lewin.
\newblock Fields of fractions for group algebras of free groups.
\newblock {\em Transactions of the American Mathematical Society},
  192:339--346, 1974.

\bibitem[LR13]{lau13}
Aaron Lauve and Christophe Reutenauer.
\newblock Rational series in the free group and the connes operator.
\newblock {\em Contemporary Mathematics}, 592:177--197, 2013.

\bibitem[McC80]{mcc80}
James McCool.
\newblock A characterization of periodic automorphisms of a free group.
\newblock {\em Transactions of the American Mathematical Society},
  260(1):309--318, 1980.

\bibitem[Mor12]{mor12}
Primo{\v{z}} Moravec.
\newblock Groups of order $p^{6}$ and their unramified brauer groups.
\newblock {\em Journal of Algebra}, 372:420--427, 2012.

\bibitem[Myr75]{myr75}
Jan Myrheim.
\newblock A theorem on restricted group representations.
\newblock {\em Mathematica Scandinavica}, 37:193--196, 1975.

\bibitem[PV94]{pop94}
Vladimir~L Popov and Ernest~B Vinberg.
\newblock Invariant theory.
\newblock In {\em Algebraic geometry IV}, pages 123--278. Springer, 1994.

\bibitem[Reu99]{reu99}
Christophe Reutenauer.
\newblock Malcev-neumann series and the free field.
\newblock {\em Expositiones Mathematicae}, 17(5):469--478, 1999.

\bibitem[Sal84]{sal84}
David~J Saltman.
\newblock Noether's problem over an algebraically closed field.
\newblock {\em Inventiones mathematicae}, 77(1):71--84, 1984.

\bibitem[Spr06]{spr06}
T.A. Springer.
\newblock {\em Invariant Theory}.
\newblock Lecture Notes in Mathematics. Springer Berlin Heidelberg, 2006.

\bibitem[Swa83]{swan83}
Richard~G Swan.
\newblock Noether’s problem in galois theory.
\newblock In {\em Emmy Noether in Bryn Mawr}, pages 21--40. Springer, 1983.

\bibitem[Vol18]{Vol18}
Jurij Vol\v{c}i\v{c}.
\newblock Matrix coefficient realization theory of noncommutative rational
  functions.
\newblock {\em J. Algebra}, 499:397--437, 2018.

\end{thebibliography}

\bibliographystyle{alpha}

\end{document}